\documentclass[oneside, 12pt]{amsart}
\usepackage{amscd, amssymb, amsmath, mathrsfs}
\usepackage[english]{babel}
\usepackage{booktabs}
\usepackage{tikz-cd}
\usepackage{url}
\usepackage[pdftex, colorlinks=true,  citecolor=blue, linkcolor=blue, linktocpage=true]{hyperref}

\newcommand\mylabel[1]{\label{#1}\marginpar{\vspace{-1ex}\medskip\medskip\footnotesize \tt #1}}
\renewcommand\mylabel[1]{\label{#1}}
\newcommand{\mydate}{
\number\day\space
\ifcase\month \or January\or February\or March\or April\or May\or June\or July\or August\or September\or October\or November\or December\fi 
\space\number\year}

\DeclareUrlCommand\arXiv{\urlstyle{same}}

\newtheorem{theorem}{Theorem}[section]
\newtheorem*{maintheorem}{Theorem}
\newtheorem{lemma}[theorem]{Lemma}
\newtheorem{proposition}[theorem]{Proposition}
\newtheorem{corollary}[theorem]{Corollary}

\theoremstyle{definition}
\newtheorem{definition}[theorem]{Definition}

\newtheorem*{acknowledgement}{Acknowledgement}

\theoremstyle{remark}

\newcommand{\ZZ}{\mathbb{Z}}

\newcommand{\PP}{\mathbb{P}}
\renewcommand{\AA}{\mathbb{A}}
\newcommand{\GG}{\mathbb{G}}

\newcommand{\ideala}{\mathfrak{a}}

\newcommand{\shA}{\mathscr{A}}

\newcommand{\shC}{\mathscr{C}}

\newcommand{\shI}{\mathscr{I}}

\newcommand{\shL}{\mathscr{L}}

\newcommand{\catC}{\mathcal{C}}

\newcommand{\foP}{\mathfrak{P}}

\newcommand{\foX}{\mathfrak{X}}

\newcommand{\aff}{\text{\rm aff}}
\newcommand{\Aff}{\text{\rm Aff}}
\newcommand{\alg}{\text{\rm alg}}

\newcommand{\Alb}{\operatorname{Alb}}

\newcommand{\Ass}{\operatorname{Ass}}
\newcommand{\Aut}{\operatorname{Aut}}

\newcommand{\Br}{\operatorname{Br}}

\newcommand{\can}{\text{\rm  can}}

\newcommand{\Cokernel}{\operatorname{Coker}}

\newcommand{\Div}{\operatorname{Div}}

\newcommand{\Exc}{\operatorname{Exc}}

\newcommand{\Fitt}{\operatorname{Fitt}}

\newcommand{\Gal}{\operatorname{Gal}}

\newcommand{\Hom}{\operatorname{Hom}}

\newcommand{\id}{{\operatorname{id}}}
\newcommand{\Image}{\operatorname{Im}}

\newcommand{\Kernel}{\operatorname{Ker}}

\newcommand{\Lie}{\operatorname{Lie}}

\newcommand{\dirlim}{\varinjlim}
\newcommand{\invlim}{\varprojlim}

\newcommand{\lra}{\longrightarrow}

\newcommand{\maxid}{\mathfrak{m}}

\newcommand{\primid}{\mathfrak{p}}
\renewcommand{\O}{\mathscr{O}}

\newcommand{\op}{\text{\rm op}}

\newcommand{\Pic}{\operatorname{Pic}}

\newcommand{\pr}{\operatorname{pr}}

\newcommand{\quadand}{\quad\text{and}\quad}

\newcommand{\ra}{\rightarrow}

\newcommand{\red}{{\operatorname{red}}}
\newcommand{\Reg}{\operatorname{Reg}}

\newcommand{\Res}{\operatorname{Res}}

\newcommand{\sep}{{\operatorname{sep}}}
\newcommand{\Set}{{\text{\rm Set}}}

\newcommand{\sh}{{\text{\rm sh}}}

\newcommand{\Sing}{\operatorname{Sing}}

\newcommand{\Spec}{\operatorname{Spec}}

\newcommand{\Supp}{\operatorname{Supp}}

\newcommand{\uHom}{\underline{\operatorname{Hom}}}

\newcommand{\lieg}{\mathfrak{g}}

\newcommand{\lieh}{\mathfrak{h}}

\newcommand{\Cpt}{\operatorname{Cpt}}
 
\newcommand{\baf}{\bar{f}}

\newcommand{\uExt}{\underline{\operatorname{Ext}}}

\begin{document}

\title[Albanese maps]
      {Albanese maps for open algebraic spaces}.

\author[Stefan Schr\"oer]{Stefan Schr\"oer}
\address{Mathematisches Institut, Heinrich-Heine-Universit\"at, 40204 D\"usseldorf, Germany}
\curraddr{}
\email{schroeer@math.uni-duesseldorf.de}

\subjclass[2010]{14K30, 14G17, 14E15, 14A20, 14L15}

\dedicatory{Second revised version, 1 July 2023}

\begin{abstract}
We show that for each algebraic space that is separated and of finite type over a field, and whose affinization
is   connected and reduced, 
there is a universal morphism to a para-abelian variety. The latter are schemes that acquire  the structure
of an  abelian variety after  some ground field extension. This generalizes  classical results   of Serre on universal morphisms
from algebraic varieties to abelian varieties. Our proof relies on   corresponding facts for
the proper case, together with the structural properties of group schemes,  removal of singularities by alterations,
and ind-objects.
It turns out that the formation of the Albanese variety commutes with base-change up to universal homeomorphisms.
We also give a detailed analysis of    Albanese maps   for   algebraic curves and algebraic groups, with special emphasis on 
imperfect ground fields.
\end{abstract}

\maketitle
\tableofcontents

\section*{Introduction}
\mylabel{Introduction}

The \emph{Albanese variety} and the \emph{Albanese map} are   fundamental objects in algebraic geometry.
Originally, the construction was  purely  transcendental, depending on  path integrals over closed holomorphic one-forms.
In the form given by Blanchard \cite{Blanchard 1956}, it is the universal holomorphic map $f:X\ra \Alb(X)$ 
from a compact connected  complex space $X$ endowed
with a base-point $x_0$ to a complex torus $A=\Alb(X)$, where the image of $x_0\in X$ is the zero element $0\in A$.

Over   arbitrary ground fields $k$, Albanese maps for proper varieties and schemes were constructed 
by Matsusaka \cite{Matsusaka 1953} and Grothendieck \cite{FGA VI},
by using the Picard scheme. 
In the absence of rational points $x_0\in X$, however,   notorious complications arise.
These problems are particularly severe over imperfect fields $k$ of characteristic $p>0$, when $X$ may become geometrically non-reduced.
From my perspective, it is important to have the Albanese map in  full generality,   over imperfect fields and for
geometrically non-reduced schemes, and without the 
troublesome burden of base-points,
to apply it  in the  theory of group schemes and their torsors, and also for questions on generic fibers in Mori fibrations,
for example related to \cite{Brion; Schroeer 2023} or \cite{Fanelli; Schroeer 2020a}.

To circumvent these issues, and to clarify the classical situation as well, one has to replace   abelian varieties by the so-called \emph{para-abelian varieties}.
The latter are schemes $P$ such that for 
some field extension $k\subset k'$, the base-change $P'=P\otimes_kk'$ admits the structure of an abelian variety.
Apparently, the name was coined by Grothendieck, but did not gain widespread use (\cite{FGA VI}, Theorem 3.3).  
Roughly speaking, these schemes are like abelian varieties, but may lack rational points and group laws.
This notion, which I find very clarifying,  was  in the above form introduced and  analyzed in \cite{Laurent; Schroeer 2021}, where 
we  constructed for any proper algebraic space $X$ with $h^0(\O_X)=1$ a universal morphism $f:X\ra\Alb_{X/k}$
to a para-abelian variety. The defining property of this \emph{Albanese map}  is that
it induces an isomorphism between \emph{maximal abelian subvarieties} inside the Picard schemes.  
These Albanese maps have excellent properties: They are    functorial in $X$, commute  with ground field extensions $k\subset k'$, and 
are equivariant with respect to the action of the group scheme $\Aut_{X/k}$.

After the completion of the present work, I learned that Albanese maps where also constructed  in the setting of 
algebraic stacks by Brochard (\cite{Brochard 2021}, Section 7 and 8). They take values in commutative group stacks
that combine abelian varieties and finite group schemes, and exist under the condition that  $\Pic^\tau_{X/k}$ is proper
(loc.\ cit. Theorem 8.1).  

The goal of this paper is to \emph{extend the theory of Albanese maps by removing the assumption that $X$ is universally closed}.
In other words, given a scheme or more generally an  algebraic space $U$ that is separated and of finite type, we seek a universal morphism to a para-abelian variety.
The first main result of this paper gives a rather comprehensive answer:

\begin{maintheorem}
(See Thm.\ \ref{existence albanese map})
Let $k$ be an arbitrary ground field, and $U$ be an algebraic space that is separated and of finite type.
Suppose the affinization   $U^\aff=\Spec\Gamma(U,\O_U)$ is connected and reduced, and   that $k$ coincides with the essential   field of constants for $U$.
Then there is a universal morphism $f:U\ra\Alb_{U/k}$ to a para-abelian variety, which is functorial in $U$.
\end{maintheorem}

This extends results of Matsusaka \cite{Matsusaka 1953} and Serre  \cite{Serre 1960}    on universal maps from algebraic varieties to abelian varieties, 
obtained  in   classical language, compare also the discussion of Esnault, Srinivas and Viehweg \cite{Esnault; Srinivas; Viehweg 1992}. It also generalizes more 
recent results of Wittenberg (\cite{Wittenberg 2008}, Appendix) and Achter, Casalaina-Martin and Vial (\cite{Achter; Casalaina-Martin; Vial 2019}, Appendix in the first arXiv version),
where geometrically reduced schemes are treated. The  latter had in the mean time also obtained further results  \cite{Achter; Casalaina-Martin; Vial 2022},  which also show that our assumption
on the affinization $U^\aff$ and $k$ are inevitable. Also note that there  are numerous generalizations in connection with Albanese maps,
for example involving a modulus (for example  Serre \cite{Serre 1975}, \"Onsiper \cite{Onsiper 1989}, Russell \cite{Russell 2013}), or in relation to  cycles 
(I just mention Samuel \cite{Samuel 1960}, Murre \cite{Murre 1985}, Bloch \cite{Bloch 1976} and Kahn \cite{Kahn 2021}), but such matters are beyond the scope of this paper.

Our assumption on the \emph{essential  fields of constants}, 
a     concept introduced in Section \ref{Compactifications} that appears to be of independent interest,
ensures that for all compactifications $U\subset X$ we have $h^0(\O_X)=1$.
Note that this automatically holds after passing to a finite field extension  inside   the ring $\Gamma(U,\O_U)$.

The  main idea  for the proof of the theorem 
is    conceptual and direct: We consider   \emph{compactifications} $i_\lambda:U\ra X_\lambda$ and the resulting \emph{maximal abelian subvarieties} $A_\lambda$
inside the  Picard schemes $\Pic_{X_\lambda/k}$. The existence of compactifications
goes back to Nagata \cite{Nagata 1962}; for algebraic spaces this was more recently established by  
Conrad, Lieblich and Olsson \cite{Conrad; Lieblich; Olsson 2012}. 
The collection  $(X_\lambda,i_\lambda)_{\lambda\in L}$ of all compactifications  is, up to isomorphism, a filtered ordered set, and  
we get an  ind-object of abelian varieties  $(A_\lambda)_{\lambda\in L}$. Note that the concept of \emph{ind-objects} plays a crucial 
role in several category-theoretic constructions of Grothendieck.
We then use results of de Jong \cite{de Jong 1997} on removal of singularities by alterations,
together with other results on the behavior of Picard schemes, 
to show that    $(A_\lambda)_{\lambda\in L}$ is \emph{essentially constant}. It follows that for sufficiently large $\lambda$,
the Albanese varieties $\Alb_{X_\lambda/k}$ do  not depend on the index.
These    give  the desired Albanese variety $\Alb_{U/k}$, whose universal property easily follows from the   theory of  rational maps. 

Recall that in differential topology, an  \emph{open manifold} is a manifold without boundary whose connected components are non-compact.
In analogy, one may call an algebraic space   that is separated and of finite type but not  proper
an ``open  algebraic space'', as occurring in the title.

There is  a crucial difference between the proper and the open situation:
In the latter case, the category $\Cpt(U)$ of all compactifications usually changes under
ground field extension $k\subset k'$ in a significant way, for lack of initial object.
It is therefore   a priori unclear how the Albanese map behaves under ground field
extension.  Our second main result clarifies this:

\begin{maintheorem}
(See Thm.\ \ref{base change})
For arbitrary $k\subset k'$, the comparison map $\Alb_{U\otimes k'/k'}\ra \Alb_{U/k}\otimes k'$ is a finite universal homeomorphism.
It is an isomorphism provided that the extension $k\subset k'$ is separable.
\end{maintheorem}

We shall see that over imperfect fields $k$,
the Albanese variety   may indeed change upon inseparable extensions. This phenomenon already appears for algebraic curves $C$,
if the \emph{canonical compactification} $\bar{C}$ by regular points at infinity is  not geometrically regular at infinity,
as we show in  Theorem \ref{base change fails}.  
 
Our theory of Albanese maps also has   consequences for \emph{algebraic groups}, that is, 
group schemes $G$ of finite type.
For smooth $G$ the existence of an Albanese map is a classical result, obtained in various degrees of generality
by Barsotti \cite{Barsotti 1955}, Rosenlicht \cite{Rosenlicht 1956} and Chevalley \cite{Chevalley 1960}.
Modern accounts are given by Conrad \cite{Conrad  2002} and Brion \cite{Brion 2017}, but the general case was apparently not covered so far.
Each  algebraic group $G$ sits in a central extension $0\ra N\ra G\ra G^\aff\ra 1$,
where the kernel $N$ of the affinization map is   \emph{anti-affine}, a notion introduced and analyzed by Brion \cite{Brion 2009}.
Our theory of Albanese varieties applies provided that 
 $G^\aff$ and equivalently $G$ are  reduced and connected.

\begin{maintheorem}
(See Thm.\ \ref{albanese for group schemes})
Suppose that $G$ is reduced and connected. Then 
$\Alb_{G/k}=N/N'$, where $N'\subset N$ is the smallest subgroup scheme such that $N/N'$ is proper and
the induced projection $G/N'\ra G^\aff$ admits a section.
\end{maintheorem}

Note that  the section, if it exists,  does not necessarily respect the group laws. So
the group scheme $G/N'$ has as underlying scheme   $N/N'\times G^\aff$, and the group law arises from the product group law
by modifying it with   a Hochschild cocycle from $Z^2(G^\aff,G/N')$.
The result seems to be relevant in connection with the \emph{pseudo-abelian varieties}. These are certain extensions
of smooth connected unipotent group schemes by abelian varieties, introduced and studied by Totaro \cite{Totaro 2013}.
Using reduced connected  unipotent group schemes $U$, supersingular abelian varieties $N$, and Hochschild cohomology
we   construct in Proposition \ref{albanese disrespects group law} algebraic groups  whose Albanese map does not respect the group law.

\medskip
The paper is structured as follows:
In Section \ref{Generalities} we recall  generalities on compactifications, Macaulayfications, and resolution of singularities by alterations,
as well as the theory of para-abelian varieties.
Section \ref{Modifications} contains an analysis of the cokernels of Picard schemes for modifications $f:X\ra Y$,
in particular if $Y$ is regular.
This is used in Section \ref{Abelian part}, to understand the effect on the  abelian part of the Picard scheme.
In Section \ref{Ind-objects} we briefly recall the notion of ind-objects, and give a characterization of
ind-objects of abelian varieties that are essentially constant. Using this, we construct in Section \ref{Compactifications}
the Albanese map for an algebraic space $U$ that is separated and of finite type, in terms of
the system of its compactifications $(X_\lambda,i_\lambda)$.
In Section \ref{Behavior} we study the behavior of the Albanese variety under ground field extensions.
The last two sections study algebraic curves and algebraic groups, respectively. 

\begin{acknowledgement}
I wish to thank Brian Conrad and the referee for careful reading and  many valuable remarks, which helped to improve the paper,
and Sandor Kov\'acs for pointing out
a gap in my application of his results in \cite{Kovacs 2017}.
This research was conducted in the framework of the   research training group
\emph{GRK 2240: Algebro-Geometric Methods in Algebra, Arithmetic and Topology}, which is funded
by the Deutsche Forschungsgemeinschaft. 
\end{acknowledgement}

\section{Recollections and generalities}
\mylabel{Generalities}

Let $S$ be a base scheme, and write $(\Aff/S)$ for the category of   $S$-schemes that are affine.
Recall that an \emph{algebraic space} is a contravariant functor $X:(\Aff/S)\ra(\Set)$
that satisfies the sheaf axiom with respect to the \'etale topology, for which 
the diagonal monomorphism $X\ra X\times X$ is relatively representable by schemes,
and such that there is an \'etale surjection $U\ra X$ from some scheme $U$.
Algebraic spaces are important generalizations of schemes that in many situations allow more freely
the formation of quotients. Those that are representable by schemes are called \emph{schematic}.
We refer to the monographs of Artin \cite{Artin 1973}, Knutson \cite{Knutson 1971},   
Olsson \cite{Olsson 2016} and the stacks project \cite{SP} for a comprehensive treatment.

Let $X$ be an algebraic space. Although  not immediate from the definition, it comes with 
a topological space, and I want to discuss this matter  first:
A \emph{point}   is an equivalence class of some morphism $a:\Spec(K)\ra X$, where $K$ is a field, and the equivalence relation is 
generated by the factorization relation.
The set of all points is denoted by  $|X|$. It is endowed with the \emph{Zariski topology}, which  is the finest topology that renders 
all maps $|U|\ra |X|$ continuous, where $U\ra X$ runs through the \'etale maps from schemes $U$. 
If  our algebraic space is  schematic,   each $a\in|X|$ has via the image $x\in X$  a canonical representation  by $\Spec\kappa(x)\ra X$,
so the above  $|X|$ can be identified with the usual underlying topological space of $X$.

A \emph{geometric point} is a morphism $\bar{a}:\Spec(\Omega)\ra X$ for some algebraically closed field $\Omega$.
Note that here we do not pass to equivalence classes.
Such  morphisms lift  through every given \'etale surjection $U\ra X$ from a scheme $U$. The lift is not unique, in general, but  gives a point $u\in U$,
and the relative separable closure for the inclusion  $\kappa(u)\subset\Omega$ depends only on the point $a\in |X|$ represented
by the geometric point $\bar{a}$. We denote this by  $\kappa(a)^\sep$. Likewise, we write $\O_{X,a}^\sh$ for the
resulting strictly henselian local ring, where the residue field    $\kappa(a)^\sep$ is separably closed.

Let $Y$ be a noetherian algebraic space. Then there is a dense open subspace  $U\subset X$ that is isomorphic to a scheme
(\cite{Olsson 2016}, Theorem 6.4.1).  As explained above,   we may regard the generic points   $\eta\in |X|$ as generic points $\eta\in X$.
A \emph{modification} is a proper morphism $f:X\ra Y$ such that $f^{-1}(V)\ra V$ is an isomorphism
for some dense open set $V\subset Y$, and that $f$ induces a bijection between the sets of generic points.
These are just the proper birational morphisms, in case $Y$ is integral.
An \emph{alteration} is a proper morphism $f:X\ra Y$ such that $f^{-1}(V)\ra V$ is a finite surjection for some dense open set $V\subset Y$,
and that $f$ induces a bijection between the sets of generic points.
Let us now recall and collect three     deep  and     fundamental results:

\begin{theorem}
\mylabel{fundamental results}
Suppose the base scheme  $S=\Spec(R)$ is the spectrum of a noetherian ring.
Let $Y$ be an algebraic space that is separated and of finite type over $S$.
\begin{enumerate}
\item There is an open  embedding of algebraic spaces $Y\subset\bar{Y}$ with $\bar{Y}$ proper.
\item If the ring $R$ admits a dualizing complex, there is a modification $f:X\ra Y$ with  $X$ Cohen--Macaulay.
\item If $R$ is excellent of dimension $\leq 2$ and $Y$ is reduced, there is an alteration $f:X\ra Y$ with $X$  regular.
\end{enumerate}
In the above, one moreover may choose $X$ with an ample invertible sheaf.
\end{theorem}

The embedding  in statement (i) is usually called \emph{Nagata compactification}. The above  general form  is due
to   Conrad, Lieblich and Olsson (\cite{Conrad; Lieblich; Olsson 2012}, Theorem 1.2.1).
Note that the case of schemes was already treated by L\"utkebohmert \cite{Luetkebohmert 1993}.
The other two statements are reduced to the case where $Y$ is a scheme with Chow's Lemma,
which was established by Rydh for algebraic spaces (\cite{Rydh 2015}, Theorem 8.8).
Morphisms $f:X\ra Y$ as in    (ii)  are called \emph{Macaulayfications}. Under certain assumptions, such maps where
first constructed by Faltings \cite{Faltings 1978}. The above general form was established by  Kawasaki \cite{Kawasaki 2000}.
Further generalizations, without   assumptions on the dualizing complex of $R$,  were recently obtained by   \v{C}esnavi\v{c}ius \cite{Cesnavicius 2018}.
Result  (iii) is due to de Jong (\cite{de Jong 1997}, Corollary 5.15). The case of  ground fields was already
established earlier (\cite{de Jong 1996}, Theorem 4.1). Throughout the paper, we will freely use the above facts.

A scheme $P$ over some field $k$ is called a \emph{para-abelian variety} if some base-change $P'=P\otimes_kk'$ admits the
structure of an abelian variety. 
This  notion seems to go back to Grothendieck (\cite{FGA VI},  Theorem 3.3), in somewhat different but equivalent form, 
and was thoroughly studied in \cite{Laurent; Schroeer 2021}. It  turns out that the para-abelian varieties $P$
are precisely the torsors over abelian varieties $G$, which is the more traditional point of view.
The crucial fact here is  this additional structure of  $G$
and its action on $P$ can be reconstructed, in an intrinsic way,   from the para-abelian variety: Inside $\Aut_{P/k}$ the group scheme $G$ is the
inertia subgroup scheme with respect to   $\Pic^\tau_{P/k}$ (loc.\ cit.\ Theorem 5.3, which actually works in a relative setting).
Note that  the para-abelian varieties $P$ are projective, smooth, and connected, but may lack group laws and rational points.
In fact,  the  group laws on base-changes $P'$ correspond to the elements  $e'\in P(k')$ (loc.\ cit.\ Proposition 4.3).
Also note that $\Pic^0_{P/k}=\Pic^\tau_{P/k}$, according to \cite{Mumford 1970}, Corollary 2 on page 178.

Let $X$ be an algebraic space over some base scheme $S$ whose structure morphism $X\ra S$ is proper, flat, of finite presentation
and cohomologically flat in degree $d=0$. The latter means that $f_*(\O_X)$ is locally free of finite rank,
and that its formation commutes with base-change.
Then the sheafification of the functor $R\mapsto\Pic(X\otimes R)$ with respect to the fppf topology
is representable by an algebraic space $\Pic_{X/S}$, which is  locally of finite presentation  (\cite{Artin 1969}, Theorem 7.3).
Moreover, the subsheaf $\Pic^\tau_{X/S}$ stemming from fiberwise numerically trivial sheaves
is representable by an algebraic space that is  of finite presentation, and the inclusion is an open embedding (\cite{Laurent; Schroeer 2021}, Theorem 2.1).
This subsheaf is stable with respect to the action of the relative automorphism group scheme.

A \emph{family of para-abelian varieties} is a proper, flat morphism  $P\ra S$  of   finite presentation, where the total space $P$
is an algebraic space, and all fibers  are para-abelian varieties. 
Then the subgroup scheme $G\subset\Aut_{P/S}$ that acts trivially on $\Pic^\tau_{P/S}$   is a 
\emph{family of abelian varieties} (also known as abelian schemes),   its action on $P$ is free and transitive, and we have an identification $\Pic^\tau_{P/S}=\Pic^\tau_{G/S}$
(\cite{Laurent; Schroeer 2021}, Section 5). 
Here  $\Aut_{P/S}$ denotes the relative automorphism group scheme for $P$ over the base scheme $S$.
A morphism $f:X\ra P$ to a family of para-abelian varieties is called an \emph{Albanese map} if the resulting
$f^*:\Pic^\tau_{P/S}\ra \Pic_{X/S}^\tau$ is a monomorphism and identifies the abelian varieties $A_s=\Pic^\tau_{P/S}\otimes\kappa(s)$
with the \emph{maximal abelian subvarieties} (\cite{Laurent; Schroeer 2021}, Section 7)
inside the group schemes  $G_s=\Pic_{X/S}^\tau\otimes\kappa(s)$, for all points  $s\in S$.
If it exists, it is universal for morphisms into families of para-abelian varieties. 
We then set $\Alb_{X/S}=P$ and call it the \emph{family of Albanese varieties}.
Over ground fields,   the  existence is automatic (\cite{Laurent; Schroeer 2021}, Corollary 10.5, compare also \cite{Brochard 2021} Theorem 8.1):

\begin{theorem}
\mylabel{albanese for proper}
If $S=\Spec(k)$ is the spectrum of a field, then every proper algebraic space $X$ with $h^0(\O_X)=1$
has an Albanese map $X\ra\Alb_{X/k}$. Moreover, the formation of the Albanese variety $\Alb_{X/k}$ is functorial
in $X$,   equivariant with respect to the action of the group scheme $\Aut_{X/k}$, and commutes with ground field extensions.
\end{theorem}

Note that indeed the group scheme $\Aut_{X/k}$, which in positive characteristics could be non-regular and even non-reduced,
and not only its group of rational points  acts on
the Albanese variety, thanks to our treatment of the relative setting (\cite{Laurent; Schroeer 2021}, Corollary 10.3).

Before we proceed, recall that a group scheme $U$  over a ground field $k$
is called \emph{unipotent} if  the base change $U\otimes k^\alg$
admits  a composition series whose quotients embed into the additive group $\GG_{a,k^\alg}$. For more details,
we refer to \cite{SGA 3b}, Expos\'e XVII, Section 1.
 
Now let  $h:P_1\ra P_2$ be a morphism between para-abelian varieties over a ground field $k$.
Write $G_i\subset\Aut_{P_i/k}$ for the subgroup schemes that act trivially on  $\Pic^\tau_{P_i/k}$.
Then  there is 
a unique homomorphism $h_*:G_1\ra G_2$ between these abelian varieties that makes $h$ equivariant with respect to the resulting $G_1$-actions 
(\cite{Laurent; Schroeer 2021}, Proposition 5.4).
Let $A_i=\Pic^\tau_{G_i/k}$ be the \emph{dual abelian varieties}, and $h^*:A_2\ra A_1$ be the induced homomorphism.
These maps are related as follows:

\begin{proposition}
\mylabel{equivalences abelian varieties}
In the above situation, the following equivalences holds: 
\begin{enumerate}
\item $h$ is surjective $\Leftrightarrow$ $h_*$ is surjective $\Leftrightarrow$ $h^*$ is finite.
\item $h$ is finite $\Leftrightarrow$ $h_*$ is finite $\Leftrightarrow$ $h^*$ is surjective.
\item $h$ has geometrically connected fibers $\Leftrightarrow$  $\Kernel(h_*)^\aff$ is local $\Leftrightarrow$ $\Kernel(h^*)^\aff$ is unipotent.
\end{enumerate}
\end{proposition}

\proof
It suffices to treat the case that $k$ is algebraically closed. Then there is a rational point $a_1\in P_1$
giving   identifications $P_i=G_i$. So $h_*=h$, and we  thus may start with a homomorphism $h:G_1\ra G_2$ of abelian varieties.
Moreover,  $\Image(h)$ and $\Cokernel(h)$ are abelian varieties,
and   $N=\Kernel(h)$ is proper. 
These are related by  three short exact sequences:
\begin{gather*}
0\lra G\lra N\lra N/G\lra 0,\\
0\lra \Image(h)\lra G_2\lra \Cokernel(h)\lra 0,\\
0\lra N\lra G_1\lra \Image(h)\lra 0.
\end{gather*}
The kernel $G$ of the affinization map $N\ra N^\aff$ is connected and smooth 
(\cite{Demazure; Gabriel 1970}, Chapter III, 8.2). It is also proper because $N$ is proper.
Thus $N$ is an extension of the finite group scheme $N/G$ by the abelian variety
$G$.  
Now recall that for any abelian variety $B$, the  dual abelian variety  $B^*=\Pic^\tau_{B/k}$ represents the sheaf
$\uExt^1(B,\GG_m)$, as explained in \cite{Moret-Bailly 1981}, Appendix, and that $\uHom(B,\GG_m)=0$.
Moreover,  for any finite group scheme $H$, the \emph{Cartier dual} $H^\vee$ represents the sheaf $\uHom(H,\GG_m)$.

Applying the contravariant functor $\uHom(\cdot,\GG_m)$ to the above short exact sequences and using Lemma \ref{sheaf ext} below, we obtain   identifications
$\uExt^1(N,\GG_m)=G^*$ and $\uHom(N,\GG_m)=(N/G)^\vee$, 
together with    exact sequences $0\leftarrow \Image(h)^*\leftarrow G_2^*\leftarrow \Cokernel(h)^*\leftarrow 0$ and 
$0\leftarrow \uExt^1(N,\GG_m)\leftarrow G_1^*\leftarrow \Image(h)^*\leftarrow\uHom(N,\GG_m)\leftarrow 0$.
In turn, we get a commutative  diagram  where the sequences with kinks are exact:
$$
\begin{tikzcd}[row sep=tiny, column sep=small ]
0	& G^*\ar[l]	& G_1^*\ar[l]		&		& G_2^*\ar[ld]\ar[ll,dashed,"h^*"']	& \Cokernel(h)^*\ar[l]	& 0\ar[l]\\
	&	&	& \Image(h)^*\ar[ld]\ar[lu]\\
	&	& 0	&		& (N/G)^\vee\ar[lu]	& 0\ar[l]
\end{tikzcd}
$$
Using the Snake Lemma, we see that $\Kernel(h^*)$ is an extension of the abelian variety $\Cokernel(h)^*$ by the finite group scheme
$(N/G)^\vee$. Consequently $h$ is surjective if and only if $\Kernel(h^*)$ is finite, which gives (i).
Assertion (ii) follows from biduality (\cite{Mumford 1970}, Corollary on page 132). Equivalently, one may argue that $h$ is finite if and only if $G^*=0$,
which means that $h^*$ is surjective.

It remains to check (iii). First observe that  $h:G_1\ra G_2$ has  geometrically connected fibers if and only if the finite group scheme $H=N/G$ is connected.
Decompose $H=H_m\oplus H_u$ into its multiplicative and unipotent   parts, and furthermore  
$$
H_m=H_m^0\oplus H_m^e\quadand H_u=H_u^0\oplus H_u^e
$$
into connected and \'etale parts. Such decompositions indeed exist  since $k$ is algebraically closed
(\cite{Demazure; Gabriel 1970}, Chapter IV, \S3, Theorem 1.1).
The Cartier duals of the connected group schemes  $H_m^0$ and $H_u^0$ are unipotent,
whereas the \'etale parts $H_m^e$ and $H_u^e$    have multiplicative Cartier duals. This gives (iii).
\qed

\medskip
In the above proof, we have used the following fact:

\begin{lemma}
\mylabel{sheaf ext}
For each inclusion of abelian varieties $A\subset B$ over a ground field $k$, the  induced  map $\uExt^1(B,\GG_m)\ra \uExt^1(A,\GG_m)$ is surjective.
For each finite group scheme $H$, the sheaf $\uExt^1(H,\GG_m)$   vanishes.
\end{lemma}

\proof
Note that for both statements one  may replace the ground field $k$ by any finite field extension.
So by  Poincare's Complete Reducibility Theorem (\cite{Mumford 1970}, Theorem 1 on page 173), we   may assume that there is another 
abelian subvariety $A'$ such  that $A\oplus A'\ra B$ is surjective, with finite kernel $N$. Choose an integer $n\geq 1$ that annihilates
the group scheme $N$. In the long exact sequence
$$
\ldots\lra \uExt^1(B,\GG_m)\lra \uExt^1(A,\GG_m)\oplus\uExt^1(A',\GG_m)\lra \uExt^1(N,\GG_m)\lra \ldots,
$$
the term on the right is annihilated by $n$, whereas the cokernel for the map on the left is an abelian variety.
It follows that the  map on the right is zero, thus the map on the left is surjective.
Actually, the group $\uExt^1(N,\GG_m)$ vanishes, as we are about to show.
 
Concerning $H$, we have to show that each    $0\ra \GG_{m,R}\ra E\ra H_R\ra 0$ over some ring $R$ 
splits  after base-change with respect to an fppf extension  $R\subset R'$.
According to \cite{Demazure; Gabriel 1970}, Chapter IV, \S3, Theorem 1.1
there is a short exact sequence $0\ra H'\ra H\ra H''\ra 0$ where $H'$ is unipotent and $H''$ is multiplicative.
The arguments for \cite{Laurent; Schroeer 2021}, Proposition 6.1 show that $\uExt^1(H'',\GG_m)=0$.
It remains to treat the case that $H$ is unipotent. In characteristic $p=0$ the finite group scheme  $H$ must be trivial, so 
we assume $p>0$. After enlarging the ground field $k$, we may assume that
$H$ admits a composition series whose quotients are isomorphic to $(\ZZ/p\ZZ)_k$ or $\alpha_p$, which reduces the problem
to these two particular cases. In both cases, the Kummer sequence yields a long exact sequence
$$
\ldots\lra\uExt^1(H,\mu_p)\lra \uExt^1(H,\GG_m)\stackrel{p}{\lra} \uExt^1(H,\GG_m)\lra\ldots,
$$
where the map on the right is zero. It thus suffices to verify that each extension $0\ra \mu_{p,R}\ra E\ra H_R\ra 0$
splits over some fppf extension $R\subset R'$.
For $H=(\ZZ/p\ZZ)_k$, this happens if we trivialize the fiber of $E\ra H_R$ over the 1-section, which is a $\mu_{p,R}$-torsor.
For $H=\alpha_p$, we pass   to the  extension $0\ra \alpha_{p,R}\ra E'\ra (\ZZ/p\ZZ)_R\ra 0$ of Cartier duals,
and argue in the same way.
\qed

\section{Modifications and regularity}
\mylabel{Modifications}

Let $k$ be a ground field of arbitrary characteristic $p\geq 0$, 
and $f:X\ra Y$ be a proper morphism between algebraic spaces that are separated and of finite type over the ground field $k$,
with $\O_Y=f_*(\O_X)$. 
By the projection formula, the induced map on Picard groups is injective, giving an inclusion $\Pic(Y)\subset \Pic(X)$.
The goal of this section is to study the resulting quotient.

Throughout, we are mainly interested in the case that $Y$ is     \emph{regular}. In other words,  the Krull dimension
of $\O_{Y,b}$ coincides with the dimension of the cotangent space  $\maxid_b/\maxid_b^2$, for all points $b\in Y$.
But note that $Y$ may fail to be geometrically regular, and it   actually may be geometrically non-reduced.
We start with a simple observation:
 
\begin{proposition}
\mylabel{quotient finitely generated}
If $f:X\ra Y$ is a modification with $Y$ regular and $X$ normal, then the group $\Pic(X)/\Pic(Y)$ is finitely generated.
\end{proposition}

\proof 
The exceptional locus $E=\Supp(\Omega^1_{X/Y})$ has codimension at least one, and its image $Z=f(E)$ has codimension at least two.
In turn, the inclusion $\Pic(Y)\subset\Pic(X)$ comes with a canonical retraction, sending an invertible sheaf $\shL$
to the bidual $f_*(\shL)^{\vee\vee}$. It follows that $\Pic(Y)$ has a canonical complement inside $\Pic(X)$,
given by the isomorphism classes of invertible sheaves of the form $\shL=\O_X(D)$,
where $D$ is a Cartier divisor supported by $E$. 

Write $\Div_E(X)$ for the group of      Cartier divisors supported on $E$. Then the image of  $\Div_E(X)\ra \Pic(X)$
is the   complement for $\Pic(Y)$.
Since $X$ is normal, the map $\Div(X)\ra Z^1(X)$ sending a Cartier divisor to the resulting Weil divisor is injective.
This gives an inclusion $\Div_E(X)\subset \bigoplus_{i=1}^r\ZZ E_i$, where $E_1,\ldots,E_r$ are the
irreducible components of codimension one contained in the exceptional locus $E$. Thus $\Div_E(X)$ is free and finitely generated, so
its image in $\Pic(X)$ is at least finitely generated.
\qed

\medskip
For $Y$ proper we now consider the group schemes $\Pic^\tau_{Y/k}$, which are of finite type, and likewise for $X$.
The formation of   $f_*(\O_X)$ commutes with flat base-change.
It follows that for each  $k$-algebra $R$, the map $\O_{Y\otimes R}\ra (f\otimes\id_R)_*(\O_{X\otimes R})$
is bijective as well, and the same holds for the multiplicative sheaf of units.
In turn, the map on Picard groups $\Pic(Y\otimes R)\ra \Pic(X\otimes R)$ 
is injective. Consequently, the map   $f^*:\Pic_{Y/k}\ra\Pic_{X/k}$ between sheafifications is a monomorphism.
It is thus a closed embedding by \cite{SGA 3a}, Expos\'e $\text{\rm VI}_B$, Corollary 1.4.2.
The short exact sequence 
$$
0\lra \Pic^\tau_{Y/k}\lra  \Pic^\tau_{X/k} \lra Q\lra 0
$$
defines another group scheme $Q$ of finite type.

\begin{theorem}
\mylabel{picard quotient affine}
If $f:X\ra Y$ is a modification with $Y$   proper and regular, then the group scheme  $Q$ defined above is affine.
\end{theorem}

\proof
It suffices to check this after some  ground field extension $k\subset k'$ for which $Y\otimes k'$ remains regular.
So from now on, we assume that $k$ is separably closed.
Seeking a contradiction, we assume that $Q$ is not affine.
Then the kernel  of the affinization map $Q\ra Q^\aff$ is non-zero. 
We now use the functorial three-step filtration $Q \supset Q_1\supset Q_2\supset Q_3$ introduced in \cite{Laurent; Schroeer 2021}, Section 7,
which in turn is based  on work of Brion (\cite{Brion 2009}, \cite{Brion 2017}).
Here  $Q_1$ is the kernel of the affinization map, and $Q_2$ is a smooth connected affine group scheme.
Moreover,  $Q_1$ is anti-affine, which means $h^0(\O_{Q_1})=1$, and the quotient $Q_1/Q_2$ is an abelian variety.
In our situation, the latter is non-trivial, because $Q_1$ is non-trivial and anti-affine, and $Q_2$ is affine.
 
We seek to relate this information on the quotient to the Picard scheme. The cartesian square
$$
\begin{CD}
P		@>>>	Q_1\\
@VVV			@VVV\\
\Pic^\tau_{X/k}	@>>>	Q
\end{CD}
$$
defines  another group scheme $P$, which is  of finite type and comes with an epimorphism $h:P\ra Q_1$.
This $P$ is not affine, because otherwise all its quotients, and in particular $Q_1$ must be affine.
As in the preceding paragraph we infer  that the anti-affine group scheme $P_1=\Kernel(P\ra P^\aff)$ is non-trivial. 
If the composition   $P_1\ra Q_1/Q_2$ vanishes, we get an induced
epimorphism $P^\aff=P/P_1\ra Q_1/Q_2$ from an affine group scheme to an abelian variety, thus $Q_1/Q_2=0$, contradiction.
Since $P_2$ is smooth, connected and affine, it belongs to the kernel of $P_1\ra Q_1/Q_2$,
and the image of the latter must be a non-zero abelian variety $B$, arising as a quotient of the abelian variety $P_1/P_2$.
Setting $A=P_1/P_2$ we get  a short exact sequence
$$
0\lra P_2\lra P_1\lra A\lra 0,
$$
where the abelian variety $A$ surjects onto $B \subset Q/Q_2$.

To proceed, we now fix some prime number $\ell$ different from the characteristic exponent of the ground field $k$.
For each commutative group scheme  $H$ of finite type,  the multiplication maps $\ell^n:H\ra H$, $n\geq 0$
induce  multiplication by $\ell^n$ on the Lie algebra $\Lie(H)$. If  $H$ is  furthermore smooth and connected, it  follows that the kernel $H[\ell^n]$ is finite and reduced,
hence $\ell^n:H\ra H$ is surjective, and thus an epimorphism. Furthermore, for each rational point $a\in H$ the fiber is an $H[\ell^n]$-torsor.
The torsor is trivial because $k=k^\sep$, so the fiber contains a rational point,
and we infer that $\ell^n:H(k)\ra H(k)$ is surjective.
Applying these facts to our group schemes $H=P_i$ and $H=A$, we see that the terms in 
$$
0\lra P_2[\ell^n]\lra P_1[\ell^n]\lra A[\ell^n]
$$
are finite and \'etale. Since $k=k^\sep$, they are actually constant. Furthermore, the map on the right is   surjective,
because $\ell^n:P_2(k)\ra P_2(k)$ is surjective.

Let $T_n$ be the group of rational points in $P_1[\ell^n]$. Their union $T=\bigcup_{n\geq 0} T_n$ 
is some $\ell$-divisible  group, whose image in $A$ is Zariski dense.
In turn, its image in $B\subset Q/Q_2$ is Zariski dense as well.
 Since $k$ is separably closed, we have $\Br(k)=0$,
and hence the canonical map $\Pic(X)\ra \Pic_{X/k}(k)$ is bijective.
We thus have an inclusion $T\subset \Pic(X)$, and this subgroup is not contained in $\Pic(Y)$,
because its image in $Q/Q_2$ is non-zero.

Suppose for the moment that $X$ is normal. By  Proposition \ref{quotient finitely generated}, the quotient $\Pic(X)/\Pic(Y)$ must be finitely generated, 
so the projection $T\ra \Pic(X)/\Pic(Y)$ factors over the torsion part, because $T$ is $\ell$-divisible.
The kernel $T'\subset T$ belongs to $\Pic(Y)$, and has finite index in $T$.
We thus conclude that $T\ra B$ factors over the finite group $T/T'\subset B$.
Thus the abelian variety $B\neq 0$ contains a finite set of rational points that is Zariski dense, a contradiction.

It remains to reduce the general situation to this special case. Let $g:X'\ra X$ be the normalization of $X_\red$
and set $f'=f\circ g$. The map $\O_Y\ra f'_*(\O_{X'})$ is bijective, by Zariski's Main Theorem, and our assertion applies
to $f':X'\ra Y$. The morphism $g:X'\ra X$ is proper and surjective, and we write $K$ for the kernel of
the induced homomorphism $g^*:\Pic^\tau_{X/k}\ra \Pic^\tau_{X'/k}$.
Applying the Snake Lemma to the commutative diagram
$$
\begin{CD}
0	@>>>	\Pic^\tau_{Y/k}	@>>>	\Pic^\tau_{X/k}	@>>> 	Q	@>>> 0\\
@.		@VVV			@VVV			@VVV\\
0	@>>>	\Pic^\tau_{Y/k}	@>>>	\Pic^\tau_{X'/k}	@>>> 	Q'	@>>> 0\\
\end{CD}
$$
we see that $Q$ is an extension of some subgroup scheme inside $Q'$ by the kernel  $K$, so we   have to verify that $K$ is affine.
If $X$ and $Y$ are schematic, this   follows from \cite{SGA 6}, Expos\'e XII, Corollary 1.5, and the
reasoning immediately extends to algebraic spaces.  
\qed

\medskip
If both $Y$ and $X$ are regular, I suspect that $Q=0$, but I do not know if this is always the case.
Although not needed in what follows, I like to state the following general fact:

\begin{proposition}
\mylabel{proper monomorphism}
Suppose $Y$ is proper and  $f:X\ra Y$ satisfies $\O_Y=f_*(\O_X)$ and $H^0(Y,R^1f_*(\O_X))=0$. Then  the group scheme  $Q$ is finite.
\end{proposition}

\proof
We may assume that $k$ is algebraically closed. Furthermore, the 
Leray--Serre spectral sequence yields the five-term exact sequence
\begin{equation}
\label{five term sequence}
0\ra H^1(\O_Y)\ra H^1(\O_X)\ra H^0(Y,R^1f_*(\O_X))\ra H^2(\O_Y)\ra H^2(\O_X),
\end{equation}
giving an identification $H^1(Y,\O_Y)=H^1(X,\O_X)$ and an inclusion $H^2(Y,\O_Y)\subset H^2(X,\O_X)$. 

Set $H=\Pic^0_{Y/k}$ and $G=\Pic^0_{X/k}$. We have an inclusion $H\subset G$, and our task is to show that $\dim(G)=\dim(H)$.
The  Lie algebras $\lieh=\Lie(H)$ and $\lieg=\Lie(G)$ are given by the   cohomology groups $H^1(Y,\O_Y)$ and $H^1(X,\O_X)$, respectively.
Hence $f:X\ra Y$ induces a bijection between the Lie algebras.
In characteristic zero,    we then have $h^1(\O_Y)=\dim(H)$, and likewise for $X$, so the result follows.

Suppose   $p>0$. Here we need additional arguments, which rely on  Mumford's  theory of Bockstein operations (\cite{Mumford 1966}, Lecture 27).
Since $k$ is perfect,  the reduced part  $G_\red\subset G$
is a  subgroup scheme, which must be smooth, and $\dim(G)$ coincides with the vector space dimension of $\Lie(G_\red)$.
   
Write $W_n=W_n(\O_Y)$ for the sheaf of \emph{Witt vectors} of length $n$. This sheaf of rings comes with
an additive map $V:W_n\ra W_n$ called \emph{Verschiebung}. The image of its $m$-fold iteration is denoted by  $V_n^m\subset W_n$.
As explained in \cite{Fanelli; Schroeer 2020b}, Section 2 the combination of the short exact sequences
$0\ra V^r_{r+1}\ra W_{r+1}\ra W_r\ra 0$ and $ 0\ra V^1_r\ra W_r\ra \O_Y\stackrel{\pr}{\ra}  0$ yields 
$W_r(k)$-linear maps
\begin{equation}
\label{bockstein}
\Image(H^{i}(W_r)\stackrel{\pr_*}{\ra} H^{i}(\O_Y))\;\stackrel{\beta_r}{\lra} \;\Cokernel(H^{i}(V^1_r)\stackrel{V^{r-1}_*}{\ra} H^{i+1}(V^r_{r+1})),
\end{equation}
where the image on the left is formed with respect to the  canonical projection $\pr_*W_r\ra W_1=\O_Y$,
and the cokernel on the right comes from  a     composite map $V^{r-1}:V^1\ra V^r$.
The above $\beta_r$ are called \emph{Bockstein operators}. The kernel of $\beta_r$ comprises those cohomology classes in $H^i(Y,\O_Y)$
that lift to $H^i(Y, W_{r+1})$, and we write $H^i(Y,\O_Y)[\beta]$ for their intersection $\bigcap_{r\geq 0}\Kernel(\beta_r)$.
This vector subspace has a  geometric meaning:
According to \cite{Mumford 1966}, Theorem on page 196 we have  $\Lie(G_\red)= H^1(Y,\O_Y)[\beta]$. 
It can also be seen as the intersection of the images for $H^1(Y,W_r)\ra H^1(Y,\O_Y)$, $r\geq 1$.

Our remaining task is to verify that the latter images 
 in the cohomology groups  $H^1(Y,\O_Y)=H^1(X,\O_X)$ are the same, whether
computed on $Y$ or $X$. For this, it suffices to check that the canonical maps
$H^1(Y,W_r)\ra H^1(X,W_r)$, $r\geq 1$ are bijective. We proceed by induction on $r\geq 1$.
The case $r=1$ is trivial. Suppose now $r\geq 2$, and that the assertion holds for $r-1$.
Consider the short exact sequence $0\ra V^{r-1}_r\ra W_r\ra W_{r-1}\ra 0$ on $Y$.
First note that $H^0(Y,W_r)\ra H^0(Y,W_{r-1})$ is surjective, because the map of set-valued sheaves
$W_r\ra W_{r-1}$ admits a section, and the same holds on $X$. In turn, we get a commutative diagram
$$
\begin{CD}
0	@>>>	H^1(Y,V_r^{r-1})	@>>>	H^1(Y,W_r)	@>>> 	H^1(Y,W_{r-1})	@>>>	H^2(Y,V_r^{r-1})\\
@.		@VVV			@VVV		@VVV			@VVV\\
0	@>>>	H^1(X,V_r^{r-1})	@>>>	H^1(X,W_r)	@>>> 	H^1(X,W_{r-1})	@>>>	H^2(X,V_r^{r-1})
\end{CD}
$$
with exact rows. 
The Verschiebung $V:W_n\ra W_n$ is additive, and its $(r-1)$-fold composition induces identifications
$\O=W_r/V_r^1\ra V_r^{r-1}/V_r^r=V_r^{r-1}$ on both $Y$ and $X$. Consequently,
the vertical map on the left is bijective, and the vertical map on the right is injective.
By induction, $H^1(Y,W_{r-1})\ra H^1(X,W_{r-1})$ is bijective, and we conclude with the Five Lemma (\cite{Mitchell 1965}, Chapter I, Proposition 21.1).
\qed

\medskip
In the first version of this paper, I claimed that $Q$ is finite   for any modification  $f:X\ra Y$ with $Y$ proper and regular, but
the arguments contained a gap. The statement holds if $R^1f_*(\O_X)=0$, by the above,
and this indeed is true if $X$ is  Macaulay and  also normal, according to the work of Kov\'acs \cite{Kovacs 2017}.
The existence of a modification $X'\ra X$ that is both Macaulay and normal seems to be an open problem.

\section{The abelian part of the Picard scheme}
\mylabel{Abelian part}

Let $k$ be a ground field of characteristic $p\geq 0$.
As explained in \cite{Laurent; Schroeer 2021}, Section 7
any group scheme $G$ of finite type contains a \emph{maximal abelian subvariety} $A\subset G$, which is functorial in $G$ and compatible with ground
field extensions $k\subset k'$.
Let  $Y$  be a proper algebraic space. Recall from Section \ref{Generalities} that  $\Pic^\tau_{Y/k}$ is a group scheme of finite type.
The following analogous notation seems   useful:

\begin{definition}
\mylabel{abelian part}
We  write $\Pic^\alpha_{Y/k}$ for the maximal abelian subvariety 
inside $\Pic^\tau_{Y/k}$, and call it the \emph{abelian part} of the Picard scheme.
\end{definition}

Each proper morphism $f:X\ra Y$ induces a homomorphism $f^*:\Pic^\alpha_{Y/k}\ra \Pic^\alpha_{X/k}$
of abelian varieties. We need the following:

\begin{proposition}
\mylabel{modification isomorphism}
Suppose $Y$ is  regular and that $f:X\ra Y$ is a modification. Then  the homomorphism
$f^*:\Pic^\alpha_{Y/k}\ra \Pic^\alpha_{X/k}$ of abelian varieties is an isomorphism.
\end{proposition}

\proof
Regard $P'=\Pic^\tau_{Y/k}$ as a subgroup scheme of $P=\Pic^\tau_{X/k}$,
giving an inclusion $\Pic^\alpha_{Y/k}\subset \Pic^\alpha_{X/k}$. According to Theorem \ref{picard quotient affine}, the quotient $P/P'$ is affine. Set 
$A=\Pic^\alpha_{X/k}$ and $A'=A\cap P'$. The subgroup scheme $A/A'\subset P/P'$ is likewise affine.
Being the quotient of an abelian variety, it has
$h^0(\O_{A/A'})=1$. Thus $A/A'$ is trivial, hence $A\subset P'$. The maximality of $\Pic^\alpha_{Y/k}$ yields 
$\Pic^\tau_{X/k}=A\subset \Pic^\alpha_{Y/k}$ inside $P'$, and the assertion follows.
\qed

\medskip
The main observation  in this section is the following boundedness result on the abelian part for modifications:

\begin{proposition}
\mylabel{dimension and order bounds}
There are constants $d\geq 0$ and $l\geq 1$ depending on our proper algebraic space $Y$ such that 
for every modification $f:X\ra Y$, the following holds:
\begin{enumerate}
\item 
The abelian variety $\Pic^\alpha_{X/k}$ has dimension $\leq d$.
\item
The kernel for the induced map $f^*:\Pic^\alpha_{Y/k}\ra\Pic^\alpha_{X/k}$ has order $\leq l$.
\end{enumerate}
\end{proposition}

\proof
The idea is to reduce the problem to the case that  $Y$ is regular.
Let $Y'$ be an alteration of $Y_\red$ such that       $Y'$ is regular, and let $X'$ be the reduction 
for $X\times_YY'$. This gives a commutative diagram
$$
\begin{CD}
X	@<<<	X_\red	@<<<	X'\\
@VfVV		@VVV		@VVV\\
Y	@<<<	Y_\red	@<<<	Y' 
\end{CD}
$$
where  the vertical maps are modifications.
According to Proposition \ref{modification isomorphism}, the induced map $\Pic^\alpha_{Y'/k}\ra \Pic^\alpha_{X'/k}$ is an isomorphism.
Since $X' \ra X$ is   surjective,   the kernel $K$
for   $\Pic_{X/k}\ra\Pic_{X'/k}$ is affine, according to \cite{SGA 6}, Expos\'e XII, Corollary 1.5. 
So the same holds for the induced map $A\ra A'$ on maximal abelian subvarieties. 
Its kernel $K\cap A$ is proper and affine, hence finite. Consequently $\dim(A)\leq \dim(A')$.
Thus the dimension $d\geq 0$ for the abelian parts $\Pic^\alpha_{X'/k}=\Pic^\alpha_{Y'/k}$ is the desired bound.
 
Now set $B=\Pic^\alpha_{Y/k}$ and $B'=\Pic^\alpha_{Y'/k}$. As in the preceding paragraph, the kernel $N$
for the induced map $B\ra B'$ is finite.  In light of the above commutative diagram
and the identification $\Pic^\alpha_{Y'/k}=\Pic^\alpha_{X'/k}$, we see that $N$ contains
the kernel for $\Pic^\alpha_{Y/k}\ra\Pic^\alpha_{X/k}$. Thus $l=|N|$ is the desired bound
for the kernel orders.
\qed

\section{Ind-objects of abelian varieties}
\mylabel{Ind-objects}

\newcommand{\Ind}{\operatorname{Ind}}
Throughout this section  $\catC$ denotes  some category. Let us  briefly discuss    the notation of \emph{ind-objects},
which was introduced by Grothendieck 
(\cite{SGA 4a}, Expos\'e 1, Section 8, compare also  \cite{Kashiwara; Schapira 2006}, Section 6.1).
These are  nothing but   covariant functors $A:L\ra \catC$ defined on some  filtered category $L$.
Recall that if  the category $L$ is just an ordered set, \emph{filtered}  
 means that for all $\lambda,\lambda'\in L$ there is some $\mu\in L$ with $\lambda,\lambda'\leq \mu$.
For simplicity one often writes an ind-object as $(A_\lambda)_{\lambda\in L}$, and 
calls $L$ the \emph{index category}. The morphisms $t_{\lambda\mu}:A_\lambda\ra A_\mu$ for $\lambda\ra \mu$
are usually  called \emph{transition maps}.

Each   ind-object defines a presheaf $``\!\dirlim\!"A_\lambda$ on the category  $\catC$, via the formula
$\Hom(X,``\!\dirlim\!"A_\lambda)=\dirlim\Hom(X,A_\lambda)$. Hence a morphism $``\!\dirlim\!"A_\lambda \ra ``\!\dirlim\!"B_\mu$
is a natural transformation
$$
\dirlim_\lambda\Hom(X,A_\lambda) \stackrel{\Phi_X}{\lra}  \dirlim_\mu\Hom(X,B_\mu).
$$
By the universal property of direct limits and the Yoneda Lemma, this is a compatible collection of morphisms 
$f_\lambda:\dirlim_\mu\Hom(A_\lambda,B_\mu)$. This shows 
\begin{equation}
\label{homs  ind objects}
\Hom( ``\!\dirlim_\lambda\!"A_\lambda, ``\!\dirlim_\mu\!"B_\mu)=\invlim_\lambda\dirlim_\mu \Hom(A_\lambda,B_\mu).
\end{equation}
The collection of all ind-objects, together with the above Hom sets, form the  category  $\Ind(\catC)$.
Each   $B\in \catC$ can be seen as an ind-object, with a singleton as index category, and this gives 
a fully faithful inclusion $\shC\subset\Ind(\shC)$.
By the above, a  morphism $(A_\lambda)_{\lambda\in L}\ra B$  is a compatible collection of morphisms $f_\lambda:A_\lambda\ra B$.
This morphism is an isomorphism if for some index $\lambda_0$,
there is a morphism $g:B\ra A_{\lambda_0}$ so that for all arrows $\lambda_0\ra\lambda$,
the composition
$$
B\stackrel{g}{\lra} A_{\lambda_0}\stackrel{t_{\lambda_0,\lambda}}{\lra} A_\lambda\stackrel{f_\lambda}{\lra} B
$$ 
coincides with the identity map for $B$. An ind-object $(A_\lambda)_{\lambda\in L}$ is called \emph{constant } if all transition
maps are isomorphisms.
More generally, it is called   \emph{essentially constant} if there is an index $\lambda_0$
such that for all   $\lambda_0\ra \lambda$ the transition maps $t_{\lambda_0,\lambda}:A_{\lambda_0}\ra A_\lambda$ are isomorphisms.
By abuse of notation we say that   $B=A_{\lambda_0}$ is the \emph{essential value}. Choosing for each $\lambda'$ some diagram $\lambda'\ra \lambda\leftarrow \lambda_0$,
the compositions $f_\lambda=t_{\lambda_0,\lambda}^{-1}\circ t_{\lambda',\lambda_0}$ define compatible morphisms
$f_\lambda:A_\lambda\ra B$ that    do not depend on  the choices,
and yield an isomorphism $f:``\!\dirlim\!"A_\lambda\ra B$. Thus each essentially constant ind-object
belongs to the essential image of   $\catC\subset\Ind(\catC)$.

We need the following criterion for abelian varieties: 

\begin{lemma}
\mylabel{essentially constant}
Let $(A_\lambda)_{\lambda\in L}$ be an ind-object of abelian varieties over some ground field $k$.
Suppose for each index $\lambda$, there are constants $d\geq 0$ and $l\geq 1$ such that for each arrow $\lambda\ra \mu$
we have
$$
\dim(A_\mu)\leq d\quadand |\Kernel(A_{\lambda}\ra A_\mu)|\leq l.
$$
Then the ind-object $(A_\lambda)_{\lambda\in L}$ is essentially constant.
\end{lemma}

\proof
We may replace $L$ by any cofinal subcategory. It thus suffices to treat the case that the filtered  category $L$
is just a directed ordered set, having some smallest element $\lambda$.
Since the dimensions of the $A_\mu$ are bounded, there is some $\lambda'\in L$
where $\dim(A_{\lambda'})$ takes the largest value. Replacing $\lambda$ by $\lambda'$, we may
assume that this already happens for $A_\lambda$.
Given $\lambda\leq \mu$, the transition map $A_\lambda\ra A_\mu$ has finite kernel,
and hence  $\dim(A_\lambda)= \dim(A_\mu)$.  Thus the dimensions  
are constant. In turn, the transition maps $A_\lambda\ra A_\mu$ are surjective. 

Set $A=A_\lambda$, and consider
the kernels $K_\mu=\Kernel(A\ra A_\mu)$. The Isomorphism Theorem
gives $A_\mu=A/K_\mu$.  The orders $l_\mu=|K_\mu|$
satisfy $l_\mu\leq l_\eta$ whenever $\mu\leq \eta$. On the other hand we have $l_\mu\leq l$.
Passing to some cofinal subset, we may assume that the orders $l_\mu$ are constant for all $\lambda<\mu$, so
the inclusions $K_\mu\subset K_\eta$ are equalities.
Hence the transition maps $A_\mu=A/K_\mu\ra  A/K_\eta=A_\eta$ are isomorphisms.
\qed

\section{Compactifications and Albanese maps}
\mylabel{Compactifications}

Throughout this section, $k$ is  a ground field of arbitrary characteristic $p\geq 0$, and $U$ be an algebraic space that is separated and of finite type over our ground
field $k$.
We also assume that the ring of global sections $H^0(U,\O_U)$ is  indecomposable and has trivial nil-radical.
Equivalently, the affine hull  $U^\aff=\Spec H^0(U,\O_U)$ is  connected and reduced.
Recall that if $U$ is proper with $h^0(\O_U)=1$, there is a universal map to a para-abelian variety
(\cite{Laurent; Schroeer 2021}, Corollary 10.5, confer also \cite{Brochard 2021}, Theorem 8.1).
The ultimate goal of this  section is to remove the  properness assumption.

Recall that a \emph{compactification}  is a pair $(X,i)$ where $X$ is a proper algebraic space,
and $i:U\ra X$ is an open embedding such that $X$ is the smallest closed subspace through which $i$ factors.
If $X$ is a scheme, this means that $U\subset X$ is dense and contains the finite set $\Ass(\O_X)$.
One also says that $U\subset X$ is \emph{schematically dense}.  We will apply the same locution 
for algebraic spaces.

The compactifications of $U$ form a non-empty category $\Cpt(U)$, where an arrow $(X,i)\ra (Y,j)$
is a morphism $f:X\ra  Y$ with $f\circ i= j$. We then say that $X$ \emph{dominates} $Y$.
The schematic density of $U$  ensures
that the Hom sets are empty or singletons, in other words,  the category $\Cpt(U)$ is equivalent to an ordered set.

Given   compactifications $(X_1,i_1)$ and $(X_2,i_2)$, the smallest closed subspace
through which the diagonal   $(i_1,i_2):U\ra X_1\times X_2$ factors defines another
compactification, and we see that the opposite category $\Cpt(V)^\op$ is filtered.
For each cofinal $L\subset\Cpt(V)^\op$, we thus get    ind-objects 
$$
(H^j(X_\lambda,\O_{X_\lambda}))_{\lambda\in L}\quadand 
(\Pic^\tau_{X_\lambda/k})_{\lambda\in L}\quadand
(\Pic^\alpha_{X_\lambda/k})_{\lambda\in L}
$$
taking respective values in finite-dimensional vector spaces, group schemes of finite type,
and abelian varieties. These ind-objects should be seen as invariants of interest for the algebraic space  $U$.
Before we proceed we have to address the problem of  so-called \emph{constant field extensions}.

\begin{lemma}
\mylabel{ground field essentially constant}
For each compactification $(Y,j)$ the finite $k$-algebra $H^0(Y,\O_Y)$ is a field. Moreover,
the ind-object of fields $(H^0(X_\lambda,\O_{X_\lambda}))_{\lambda\in L}$  is essentially constant.
\end{lemma}

\proof
The morphism $j:U\ra Y$ induces an inclusion $H^0(Y,\O_Y)\subset H^0(U,\O_U)$. 
Since $Y$ is proper, the $k$-algebra $F=H^0(Y,\O_Y)$ is finite. 
Hence each point in $\Spec(F)$ is generic, so the dominant map $Y\ra \Spec(F)$ is surjective.
Since the ring $H^0(U,\O_U)$ is reduced, the same holds for the subring $F$.
Using that $U$ is connected, and also dense in $Y$, we infer that $Y$ and hence its image $\Spec(F)$ is connected.
This proves the first assertion. 

Before we come to the second assertion, we make a little observation:
Let $Y'\ra Y$ be the normalization of  $Y_\red$. Consider the finite $k$-algebra  $H^0(Y',\O_{Y'})$, which
is a finite product of   finite  field  extensions of $k$.
Let $f:(X,i)\ra (Y,j)$ be a morphism from another compactification, and set $X'=(Y'\times_YX)_\red$.
We obtain a commutative diagram
$$
\begin{CD}
X'	@>>>	 X\\
@VVV		@VVfV\\
Y'	@>>>	Y.
\end{CD}
$$
The first projection   $X'\ra Y'$ is a modification of the normal algebraic space $Y'$,
and Zariski's Main theorem implies that $H^0(Y',\O_{Y'})=H^0(X',\O_{X'})$.
The second projection $X'\ra X$ induces a homomorphism $H^0(X,\O_X)\ra H^0(X',\O_{X'})$, which must
be injective because $F=H^0(X,\O_X)$ is a field. It follows that $h^0(\O_X)\leq h^0(\O_{Y'})$.
Summing up, the integers $h^0(\O_X)$ are bounded above by some number that depends only on $Y$.

This easily gives the second assertion: By passing to a cofinal set, we may assume that $L$ has a smallest member $(Y,j)$.
According to  the preceding paragraph, the numbers $h^0(\O_{X_\lambda})$ are bounded.
It follows that the ind-object  $(H^0(X_\lambda,\O_{X_\lambda}))_{\lambda\in L}$ of fields is essentially constant.
\qed

\medskip
Let us call $k'=``\!\dirlim\!" H^0(X_\lambda,\O_{X_\lambda})$  \emph{the  essential field of constants} for the algebraic space $U$.
The morphism $U\subset X_\lambda\ra (X_\lambda)^\aff=\Spec(k')$, with sufficiently large index  $\lambda$,
endows the algebraic space $U$ over $k$ with a canonical  $k'$-structure.
After replacing  $k$ by $k'$, the ground field and the essential field of constants coincide.
In what follows, we usually make this assumption, because the theory of Albanese maps for proper algebraic spaces
was developed in \cite{Laurent; Schroeer 2021}   under this condition.

\begin{proposition}
\mylabel{pic alpha essentially constant}
Suppose the ground field $k$ equals the essential field of constants for $U$.
Then the  ind-object $(\Pic^\alpha_{X_\lambda/k})_{\lambda\in L}$ of abelian varieties is essentially constant.
\end{proposition} 

\proof
By passing to a cofinal subset, we may assume that $L$ contains a smallest member $(Y,i)$.
Combining Proposition \ref{dimension and order bounds} and Lemma \ref{essentially constant}, we see that our ind-object is essentially constant.
\qed

\begin{corollary}
Suppose the ground field $k$ equals the essential field of constants for $U$.
Then there is an index $\lambda\in L$ such that for all $\mu\geq \lambda$ the transition map  
$f:X_\mu\ra X_\lambda$ induces an isomorphism $f_*:\Alb_{X_\mu/k}\ra \Alb_{X_\lambda/k}$ of para-abelian varieties.
\end{corollary}

\proof
In light of the Proposition, we may pass to  a cofinal index set and  assume that the homomorphisms $f^*:\Pic^\alpha_{X_\lambda/k}\ra \Pic^\alpha_{X_\mu/k}$ are
isomorphisms, for all $\mu\geq \lambda$ in $L$.
Let $P_\lambda=\Alb_{X_\lambda/k}$ be the Albanese varieties, and $g_\lambda:X_\lambda\ra P_\lambda$ be the Albanese maps.
By our definition of Albanese maps (\cite{Laurent; Schroeer 2021}, Section 8), 
the induced map $g_\lambda^*:\Pic^\alpha_{P_\lambda/k}\ra \Pic_{X_\lambda/k}^\alpha$ of abelian varieties
is an isomorphism.
It follows that $f_*:P_\mu\ra P_\lambda$ induces an isomorphism $\Pic^\tau_{P_\lambda/k}\ra \Pic^\tau_{P_\mu/k}$.
According to Proposition \ref{equivalences abelian varieties}, the former is an isomorphism as well.
\qed

\medskip
We now come to the main result of the paper:

\begin{theorem}
\mylabel{existence albanese map}
Suppose the ground field $k$ equals the essential field of constants for $U$.
Then there is a para-abelian variety $P$ and a morphism $f:U\ra P$ such that for each other
para-abelian variety $Q$ with a morphism $g:U\ra Q$, there is a unique morphism $h:P\ra Q$ such
that $g=h\circ f$.
\end{theorem} 

\proof
Choose a cofinal index set $L\subset\Cpt(V)^\op$ such that $f_*:\Alb_{X_\mu/k}\ra \Alb_{X_\lambda/k}$
are isomorphisms for all transition maps $f:X_\mu\ra X_\lambda$ with $\lambda,\mu\in L$, and that there is a smallest member $(Y,j)$.
Let  $P=\Alb_{Y/k}$ be the Albanese variety of the proper algebraic space $Y$.
We obtain a morphism $f:U\ra P$ as the  composition of the Albanese map $a:Y\ra P$ with the inclusion $j:U\ra Y$.

We have  to verify the universal property. Let $g:U\ra Q$ be a morphism to some other para-abelian variety.
It can be seen as a rational map $g:Y\dashrightarrow Q$.
Taking the closure of the graph $\Gamma_g\subset Y\times Q$ we obtain another compactification $(X,i)$.
The projection $\pr:X\ra Y$ is a  morphism of compactifications, and  the rational map $g:Y\dashrightarrow Q$ extends to a morphism $\tilde{g}:X\ra Q$.
We have the following commutative diagram:
$$
\begin{tikzcd} 
U\ar[r,"j"]\ar[ddr,"g"']\ar[rr, bend left=40,"i"]	& Y\ar[d,"a"]		&X\ar[l,"\pr"']\ar[ddl,"\tilde{g}"]\\
					& P\ar[d,dashed,"h"]\\
					& Q
\end{tikzcd}
$$
By construction, $a:Y\ra P$ is the Albanese map, and $\pr:X\ra Y$ induces an isomorphism between    Albanese varieties.
Consequently, the composition $a\circ \pr:X\ra P$ is the Albanese map for the proper algebraic space $X$, hence there is a unique morphism $h:P\ra Q$ 
with $h\circ a\circ\pr=\tilde{g}$. We then also have 
$h\circ a\circ j = h\circ a\circ\pr\circ i = \tilde{g}\circ i = g$,
so the whole diagram is commutative.

It remains to verify   uniqueness:  Suppose there is another morphism $h':P\ra Q$ with $h'\circ a\circ j=g$. Again the  corresponding diagram 
with $h'$ instead of $h$ is commutative,
and $h\circ a\circ j = h' \circ a\circ j$.
The compactification $j:U\ra Y$ is an epimorphism by schematic density. To see this, choose an \'etale surjection $\tilde{U}\ra U$ from 
some scheme $U'$, and apply \cite{Romagny; Rydh; Zalamansky 2018}, Lemma 2.1.1 to the resulting morphisms of schemes $\tilde{U}\ra Y$.
It follows  $h\circ a = h'\circ a$,
and in particular $h\circ a \circ \pr= h'\circ a\circ \pr$. 
Now recall that $a\circ \pr:X\ra P$ is the Albanese map for $X$. Its  universal property
ensures $h=h'$.
\qed

\medskip
By the Yoneda Lemma, the pair $(P,f)$ is unique up to unique isomorphism.
We then write $P=\Alb_{U/k}$ and call it the \emph{Albanese variety} of the algebraic space $U$.
Moreover, the morphism $f:X\ra\Alb_{U/k}$ is called the \emph{Albanese map}.
Note that if $U$ is already proper, the category $\Cpt(U)$ is equivalent to a singleton, 
hence our new construction for algebraic spaces that are separated and of finite type coincides with the
old construction for proper algebraic spaces.

\begin{proposition}
Let $g:U'\ra U$ be a morphism between   algebraic spaces that are separated and of finite type over the ground field $k$.
Suppose $k$ coincides with the essential field of constants  for both $U$ and $U'$, and that their
affine hulls are connected and reduced.
Then there is a unique morphism $g_*:\Alb_{U'/k}\ra\Alb_{U/k}$ making the diagram
$$
\begin{CD}
U'	@>g>>	U\\
@Vf'VV		@VVfV\\
\Alb_{U'/k}	@>>g_*>	\Alb_{U/k}
\end{CD}
$$
commutative, where the vertical arrows are the Albanese maps.
\end{proposition}

\proof
This is an immediate consequence of the universal property for Albanese maps. Note that the assumptions are made to 
ensure the existence of the Albanese varieties for $U$ and $U'$.
\qed

\medskip
In particular, the action of the automorphism group $\Aut(U)$   induces an action on $\Alb_{U/k}$ such that
the Albanese map is equivariant. It would be interesting to understand to what extend this holds true
for group scheme actions. Let us close the section with the following observation:

\begin{proposition}
Suppose that $U$ is connected and reduced, and that all  its irreducible components are geometrically integral.
Then the ground field $k$ coincides with the essential   field of constants for $U$.
\end{proposition}

\proof
Choose a sequence of irreducible components $U_i\subset U$, $1\leq i\leq r$  so that the successive intersections $U_i\cap U_{i+1}$ 
are non-empty, and $U=U_1\cup\ldots\cup U_r$, with repetitions allowed. Then the base-changes $U_i\otimes k^\alg$
remain   integral, and it follows that $U$ is geometrically connected and geometrically reduced.
So the  $k$-algebra $R=H^0(U,\O_U)$ is geometrically indecomposable and geometrically reduced.

We proceed by showing that $k\subset R$ is integrally closed.
Suppose there is some intermediate field $k\subset k'\subset R$.  Then $R\otimes k^\alg$ contains $k'\otimes k^\alg$. 
If $[k':k]>1$, the ring $k'\otimes k^\alg$ contains  idempotent elements $e\neq 0,1$ or     nilpotent elements $f\neq 0$, so the same
holds for the over-ring $R\otimes k^\alg$, contradiction.
Thus $k$ is integrally closed in $R$, and thus must coincide with the essential ground field for $U$.
\qed

\section{Behavior under base change}
\mylabel{Behavior}

Let $k_0$ be a ground field of characteristic $p\geq 0$, and $U_0$ be an algebraic space that is separated and of finite type over $k_0$.
Given a field extension $k_0\subset k$, we consider the base-change $U=U_0\otimes  k$.
Suppose that $U^\aff$ is reduced and connected, and that $k$ is the essential field of constants for $U$.
Then the same properties hold  for $U_0$, and we have   Albanese maps 
$$
f_0:U_0\lra\Alb_{U_0/k_0}\quadand f:U\lra \Alb_{U/k},
$$ 
and also the  base-change of the Albanese map $f_{0,k}:U\ra \Alb_{U_0/k_0}\otimes k$.
The universal property of $f$ gives a   \emph{comparison map}
\begin{equation}
\label{comparison map}
c:\Alb_{U/k}\lra \Alb_{U_0/k_0}\otimes k,
\end{equation}
such that $c\circ f =f_{0,k}$. This is an isomorphism, provided $U_0=X_0$ is proper, according to
\cite{Laurent; Schroeer 2021}, Corollary 10.5. In general, the situation is more complicated,
because the base-change functor $\Cpt(U_0)\ra \Cpt(U)$ usually is not  an equivalence of categories. 
In fact,  we shall see in the next section  examples of algebraic curves over imperfect   fields where the comparison map fails to be 
an isomorphism.  

In this section we want to establish a positive result. Recall that the field extension $k_0\subset k$ is called
\emph{separable} if for each reduced $k_0$-algebra $A_0$, the base-change $A=A_0\otimes k$ remains reduced.  

\begin{theorem}
\mylabel{base change}
In the above setting, the comparison map \eqref{comparison map} is a finite universal homeomorphism.
It is actually an isomorphism provided the field extension $k_0\subset k$ is separable.
\end{theorem}

Before entering the  proof, let us simplify notation and examine the assertions from various angles.
Set 
$$
P_0=\Alb_{U_0/k}\quadand P=\Alb_{U/k} \quadand  P_{0,k}=\Alb_{U_0/k_0}\otimes k.
$$
So our  Albanese maps are $f_0:U_0\ra P_0$ and $f:U\ra P$, and the   comparison map becomes $c:P\ra P_{0,k}$.  
Let $G\subset\Aut_{P/k}$ be the subgroup scheme that acts trivially on  $\Pic^\tau_{P/k}$.
Then $G$ is an abelian variety, its action on $P$ is free and transitive, and we have an identification $\Pic^\tau_{G/k}=\Pic^\tau_{P/k}$,
according to \cite{Laurent; Schroeer 2021}, Section 5. Similarly, we form $G_0\subset \Pic_{P_0/k_0}$ and its base-change $G_{0,k}=G_0\otimes k$. Then
there is a unique homomorphism $c_*:G\ra G_{0,k}$ making the comparison map $c:P\ra P_{0,k}$ equivariant, by loc.\ cit.\  Proposition 5.4.
We see that the assertion holds for $c$   if and only if the corresponding statement holds for $c_*$.
The latter respects the group laws,  so the assertion of the theorem means that  $c_*$ is surjective and  has local kernel.
Furthermore, $c $ and $c_*$  induce the same homomorphism  
$$
\Pic^\tau_{G_0/k_0}\otimes k=\Pic^\tau_{P_0/k_0}\otimes k
\stackrel{c^*}{\lra}
\Pic^\tau_{P/k}=\Pic^\tau_{G/k}.
$$
In light of Proposition \ref{equivalences abelian varieties}, the assertion of the theorem means that  $c^*$ is surjective and has unipotent kernel.

Suppose we have a compactification   $i_0:U_0\ra X_0$ over $k_0$, and another compactification  $i:U\ra X$ over $k$ that
dominates the base-change of $i_0$. We then get a commutative diagram 
\begin{equation}
\label{diagram compactification}
\begin{tikzcd} 
U\ar[r,"i"]\ar[rd,"i_{0,k}"']	& X\ar[r,"\baf"]\ar[d,"g"']		& P\ar[d,"c"]\\
			& X_{0,k}\ar[r,"\baf_{0,k}"']	& P_{0,k}
\end{tikzcd}
\end{equation}
of $k$-algebraic spaces,  where   $\baf_0:X_0\ra P_0$ is the Albanese map for the proper algebraic space $X_0$, such that
$f_0=\baf_0\circ i_0$, and likewise for  $\baf:X\ra P$.
The vertical map to the right is the comparison map.
By definition of Albanese maps in the proper case, the pull-back maps $\baf^*_0$ and $\baf^*$ give identifications 
$\Pic^\alpha_{X_0/k_0}=\Pic^\alpha_{P_0/k_0}$ and $ \Pic^\alpha_{X/k}=\Pic^\alpha_{P/K}$. So the assertion of the theorem means that
\begin{equation}
\label{induced map on picard}
g^*:\Pic^\alpha_{X_0/k_0}\otimes k\lra \Pic^\alpha_{X/k}
\end{equation}
is  surjective  with unipotent kernel.

\medskip
\emph{Proof of Theorem \ref{base change}.} We proceed in six   steps. Note that we may replace $X_0$ by any other compactification  $\tilde{X}_0$ of $U_0$
that dominates $X_0$, and simultaneously replace $X$ by some $\tilde{X}$ that dominates
the schematic closure of $U$ inside the fiber product $X\times_{X_0}\tilde{X}_0$.
We call this process \emph{passing to dominating compactifications}, and do this several times throughout to improve the situation.

\medskip
{\bf Step 1:}
\emph{The comparison map $c:P\ra P_{0,k}$ is surjective.} In light of Lemma \ref{equivalences abelian varieties}, this equivalently means that
the kernel of \eqref{induced map on picard}  is finite.
By construction, $g:X\ra X_{0,k}$ is surjective.
According to \cite{SGA 6}, Expos\'e XII, Corollary 1.5 the induced map  on Picard schemes has affine kernel.
Its intersection with the maximal abelian subvariety is    proper.
In turn, $\Kernel(g^*)$ is finite.  

\medskip
{\bf Step 2:}
\emph{Reduction   to the case that the field extension $k_0\subset k$ is finitely generated}.
Choose an intermediate field $k_0\subset k_1\subset k$ that is finitely generated over $k_0$,
and such that there is a proper  algebraic space $X_1$ and a  para-abelian variety $P_1$ over $k_1$, together with 
 morphisms $i_1:U_1\ra X_1$ and $\baf_1:X_1\ra P_1$
inducing the upper row in the diagram \eqref{diagram 1}, where we set $U_1=U_0\otimes k_1$.
By enlarging $k_1$ if necessary, we may assume that there are  morphisms that make the diagram 
\begin{equation}
\label{diagram 1}
\begin{CD}
U	@>i>>	X	@>\bar{f}>>		P\\
@VVV		@VVV			@VVV\\
U_1	@>i_1>>	X_1	@>\bar{f}_1>>	P_1\\
@VVV		@VVV			@VVV\\
U_0	@>>i_0>	X_0	@>>\bar{f}_0>	P_0
\end{CD}
\end{equation}
commutative.   According to \cite{Laurent; Schroeer 2021}, Proposition 8.2 the morphism 
$\baf_1:X_1\ra P_1$ is an Albanese map for the proper $k_1$-algebraic space $X_1$. 
We claim that the composite map $f_1:U_1\ra P_1$ is an Albanese map for the   $k_1$-algebraic space $U_1$.  It suffices to check that
each morphism $h_1:U_1\ra Q_1$ to some para-abelian variety, viewed as a rational map $X_1\dashrightarrow Q_1$,
is defined everywhere. In other words, the schematic closure of the  graph $\Gamma_{h_1}$ inside $X_1\times Q_1$
remains a graph. 
 By construction, this holds after base-changing along $k_1\subset k$. Since the formation of schematic closure
commutes with flat base-change, we     infer that $X_1\dashrightarrow Q_1$ is defined everywhere.

\medskip
{\bf Step 3:}
\emph{The case that the extension $k_0\subset k$ is  finite and  separable.}
This   is well-known, and we give the  arguments   for the sake of completeness:
One uses the universal property of Albanese maps to  deduce that the comparison map \eqref{comparison map} is an isomorphism.
There is a  finite extension $k\subset k'$ such that $k'$ is  Galois over both $k_0$ and $k$.
Thus it suffices to treat the case that $k_0\subset k$ is finite and  Galois. Write  $\Gamma=\Gal(k/k_0)$ for     the  Galois group, and let
$s:P\ra\Spec(k)$ be the structure morphism.
Fix some $\sigma\in \Gamma$, and consider the commutative diagram
$$
\begin{tikzcd} 
U_0\otimes k\ar[r,"\id_{U_0}\otimes \sigma"]\ar[d,"f"']\ar[dr]	& U_0\otimes k\ar[d,"f"]\\
P\ar[d,"s"']\ar[r,dashed,"\psi_\sigma"]		&  P\ar[d,"s"]\ar[dl]\\
\Spec(k)\ar[r,"\Spec(\sigma)"']		& \Spec(k).
\end{tikzcd}
$$
We now observe that the upper diagonal arrow $f\circ (\id_{U_0}\otimes\sigma)$ is a $k$-morphism
to the para-abelian variety $P$, provided that the latter is endowed with     the lower diagonal arrow $\Spec(\sigma)^{-1}\circ s$ as new structure morphism.
By the universal property of the Albanese map $f$, there is a unique dashed arrow
$\psi_\sigma:P\ra P$ making the triangles on the left commutative. 
It follows that the whole diagram is commutative. The uniqueness of $\psi_\sigma$ ensures that
the map $\Gamma\ra\Aut_{k_0}(P)$ given by $\sigma\mapsto \psi_\sigma$ respects the group laws,
and that the structure morphism $s:U\ra\Spec(k)$ is equivariant.
Since $P$ is a projective $k_0$-scheme,   the quotient $P/\Gamma$ exists as a projective $k_0$-scheme. 
We then have a canonical identification $P=(P/\Gamma)\otimes_{k_0}k$.
In particular, $P/\Gamma$ is a para-abelian variety over $k_0$.
By the above commutative diagram, the Albanese map $f:U\ra P$ over $k$ descends to a $k_0$-morphism $ U_0\ra P/\Gamma$.
Arguing as above, one sees that the latter has the universal property of the Albanese map, and infer that the comparison map must be an isomorphism.

\medskip
{\bf Step 4:}
\emph{Reduction to the cases $k=k_0(t)$ and $k=k_0(\lambda^{1/p})$.}
In light of step 2, it suffices to treat the case that $k_0\subset k$ is finitely generated.
Choose a transcendence basis $t_1,\ldots,t_n$.  Set $k_1=k_0(t_1,\ldots,t_n)$ and let $k_2$ be its relative separable closure in $k$.
The extension $k_2\subset k$ is an equality in characteristic zero. For $p>0$, it can be written
as the successive adjunction 
of certain elements $\lambda_1^{1/p},\ldots,\lambda_m^{1/p}$, according to \cite{A 4-7}, Chapter V, \S7, No.\ 7, Proposition 13.
Using inductions on $n\geq0$ and $m\geq 0$, together with step 3, we get the desired reduction.

\medskip
{\bf Step 5:}
\emph{The case $k=k_0(t)$.}
This     involves the passage to a relative setting.
For the sake of exposition, we write $F$ for the field $k=k_0(t)$, and regard it as the function field of the affine line $\AA^1_{k_0}$.
Choose a localization $k_0[t]\subset R$ by a   non-zero polynomial so that $X$
extends to a proper   morphism $s:\foX\ra\Spec(R)$. Localizing further, we may assume that $s$ is flat,
and also cohomologically flat in degree zero. In turn,   the numerically trivial part  
 $\Pic^\tau_{\foX/R}$ exists (\cite{Laurent; Schroeer 2021}, Theorem 2.1).
Doing another   localization, we may assume that $P$ extends 
to a family of para-abelian varieties $\mathfrak{P}\ra\Spec(R)$, and that the morphism $\baf:X\ra P$ extends to a morphism
$\baf_R:\foX\ra\foP$. Localizing further, we may assume that the diagram \eqref{diagram compactification} spreads out to a commutative diagram
\begin{equation}
\label{diagram 2}
\begin{tikzcd} 
U_0\otimes R\ar[r,"i_R"]\ar[rd,"i_{0,R}"']	& \foX\ar[r,"\baf_R"]\ar[d,"g_R"']		& \foP\ar[d,"c_R"]\\
			& X_0\otimes R\ar[r,"\baf_{0,R}"']	& P_0\otimes R,
\end{tikzcd}
\end{equation}
where all tensor products are over $k_0$. 
Moreover, we may assume that for each prime ideal $\primid\subset R$, the fiberwise morphisms  $U_0\otimes_{k_0}\kappa(\primid)\ra\foX\otimes_R\kappa(\primid)$
are compactifications that dominate $U_0\otimes_{k_0}\kappa(\primid)\ra X_0\otimes_{k_0}\kappa(\primid)$.
Making a final localization, we can achieve that the cokernel of $c_R$ is a family of abelian varieties,
and that the kernel is an extension of a family of finite group schemes by a family of abelian varieties.

Now recall that $R$ is a localization of the polynomial ring. Write $S=\Spec(R)$.
Then the closed points $\sigma\in S$ whose residue field $\kappa=\kappa(\sigma)$ is separable over $k_0$ form
a Zariski dense set. For such points, the morphisms $c_\kappa:\foP\otimes_R\kappa\ra P_0\otimes \kappa$ are isomorphisms by step 3.
Using flatness, we infer that $c_R$ is an isomorphism, and in particular $c=c_F$ is an isomorphism.

\medskip
{\bf Step 6:}
\emph{The case $k=k_0(\lambda^{1/p})$ in characteristic $p>0$.}
This is the most interesting and most challenging part.  Recall that by step 1 we already know that   
the comparison map $c:P\ra P_{0,k}$ is surjective. 
It remains to check that it is   injective.

First, consider the finitely many codimension-one points 
 $\zeta_1,\ldots,\zeta_r\in X_0$ 
that do not belong to $U_0$.  Passing to dominating compactifications, stemming from the semi-normalizations
of the local rings $\O_{X_0,\zeta_i}$, we may  assume that the local rings $\O_{X_0,\zeta_i}$
and $\O_{X,\zeta_i}$ are \emph{unibranch}. In other words, their henselizations have integral reductions.
Here we regard $\zeta_i$ as points in both $X_0$ and $X$, which indeed have
the same underlying topological space. 
Since $k_0\subset k$ is purely inseparable, it 
 follows that the modification $g:X\ra X_{0,k}$ is a universal homeomorphism over some open set $V_0\subset X_0$
that contains $U_0\cup\{\zeta_1,\ldots,\zeta_r\}$. 
Making a further passage to dominating compactifications, we may assume that  
$P_0=\Alb_{U_0/k_0}=\Alb_{V_0/k_0}=\Alb_{X_0/k_0}$.

Next, consider the sheaf of $\O_X$-algebras $\shA$ given by $\Gamma(W_0,\shA)=\Gamma(W_0\cap V_0,\O_{X_0})$.
This is a coherent sheaf, by \cite{EGA IVb}, Proposition 5.11.1 
and coincides with $\O_{X_0}$ on $V_0$. Now pass to dominating compactifications, where  $X_0$ is replaced by the
relative spectrum of $\shA$. It then follows that the local rings $R=\O_{X_0,a}$ satisfy Serre's condition $(S_2)$,
at each boundary point $a\in X_0\smallsetminus U_0$. In other words, the local cohomology group $H^i_{\maxid}(R)$
vanishes for $i\leq 1$. Thus the restriction map $k=H^0(X_0,\O_{X_0})\ra H^0(V_0,\O_{V_0})$ is bijective.

We now consider the finite universal homeomorphisms $g:g^{-1}(V_{0,k})\ra V_{0,k}$.
According to \cite{Kollar 1997}, Proposition 6.6 there is an integer $n\geq 0$
so that the iterated relative Frobenius map $F^n: g^{-1}(V_{0,k})\ra g^{-1}(V_{0,k})^{(p^n)}$ factors over
$V_{0,k}$. Consider the kernel $G'=G[F^n]$ for the corresponding Frobenius map $F^n:G\ra G^{(p^n)}$.
The quotient $P'=P/G'$ is para-abelian, being a  torsor for the abelian variety $G/G'$.
We claim that the diagram
$$
\begin{tikzcd} 
g^{-1}(V_{0,k})\ar[r,"f"]\ar[d, "g"']	& P\ar[d,"q"]\\
V_{0,k}\ar[r,dashed,"f'"']		& P',
\end{tikzcd}
$$
can be completed by a dashed arrow $f'$. Since the vertical arrows are   homeomorphisms, $f'$ clearly exists
as a continuous map. By the very definition of morphisms of locally ringed spaces,
we   have to check, for a given point $a\in g^{-1}(V_{0,k})$ with images  $b=f(a)$ and $b'=q(b)$,
that the composite map $\varphi:\O_{P',b'}\ra \O_{g^{-1}(V_{0,k}),a}$ factors over the subring $\O_{(V_{0,k}),g(a)}$.
The elements of  $\O_{P',b'}=k\cdot \O_{P,b}^{p^n}$ are sums of   $\lambda\cdot h^{p^n}$,
with $\lambda\in k$ and $h\in \O_{P,b}$. By our choice of $n\geq 0$, the image  $\varphi(\lambda\cdot h^{p^n})=\lambda\cdot \varphi(h)^{p^n}$
belongs to $\O_{(V_{0,k}),g(a)}$. Thus the desired $f':V_{0,k}\ra P'$ exists.

To proceed we use the existence of   \emph{Weil restrictions} along  $\Spec(k)\ra\Spec(k_0)$.
Indeed, for each $k$-scheme  $Y$ of finite type, the functor 
$$
(\Aff/k)\lra (\Set),\quad R_0\longmapsto Y(R_0\otimes_{k_0}k)
$$
is representable by 
a $k_0$-scheme $\Res_{k/k_0}(Y)$ of finite type. We refer to the monograph of Conrad, Gabber and Prasad \cite{Conrad; Gabber; Prasad 2010}, Appendix A.5
for a comprehensive treatment of Weil restrictions.
The functor respects products and closed embeddings, and therefore also group structures, torsor structures with respect to smooth group schemes, and 
being separated.
In particular $\Res_{k/k_0}(P')$ is a torsor with respect to the group scheme $\Res_{k/k_0}(G')$.
Moreover, the morphism $f':V_0\otimes k\ra P'$ corresponds to a  morphism $f'_0:V_0\ra \Res_{k/k_0}(P')$.

We already saw that $h^0(\O_{V_0})=1$, whence the image of $V_0$ in the affine hull of the Weil restriction is a rational point.
According to Proposition \ref{weil restriction} below,  the kernel of the affinization map of the group scheme $\Res_{k/k_0}(G')$ is an abelian variety. 
It follows that $f'_0$ factors over some para-abelian variety inside $\Res_{k/k_0}(P')$.
It thus also uniquely factors over a  morphism $P_0=\Alb_{V_0/k_0}\ra \Res_{k/k_0}(P')$, by the universal property of Albanese maps.
In turn,  $g^{-1}(V_{0,k})\ra P'$ factors over
some morphism $P_0\otimes k\ra P'$.
Passing to some dominating compactifications, we may assume that the composition $X\ra P\ra P'$ factors over $\Alb_{V_0/k_0}\otimes k$.
Summing up, we have a commutative diagram of para-abelian varieties
$$
\begin{tikzcd} 
	&P\ar[dl,"c"']\ar[d,"\can"]\\
P_{0,k}\ar[r]	& P'.
\end{tikzcd}
$$
Recall that the vertical arrow $P\ra P'$ is the quotient map with respect to the infinitesimal group scheme $G'=G[F^n]$,
and therefore a bijection for the underlying topological spaces.
It follows that the comparison map   $c:P\ra P_{0,k}$ is injective.
\qed

\medskip
In the above proof we have used the following fact:

\begin{proposition}
\mylabel{weil restriction}
Suppose $p>0$.
Let $k_0\subsetneqq k$ be a purely inseparable finite field extension, 
$A\neq 0$ be an abelian variety over $k$, and $G_0=\Res_{k/k_0}(A)$ the Weil restriction.
Then the affinization  $G_0^\aff$ is a smooth connected unipotent group scheme of dimension $n\geq 1$,
and the kernel of  $G_0\ra G_0^\aff$ is an abelian variety of dimension $g=\dim(A)$.
\end{proposition}

\proof
By descent, it suffices to check the properties for  $E= G_0\otimes_{k_0}k$.
According to \cite{Conrad; Gabber; Prasad 2010}, Proposition A.5.11 the  canonical homomorphism $f:E=\Res_{k/k_0}(A)\otimes_{k_0}k\ra A$ is smooth  and surjective, with
geometrically connected fibers, and $U=\Kernel(f)$ is unipotent and non-zero.
It follows that  $E$ and whence also its affinization are smooth and connected. 
Write $h:E\ra E^\aff$ for the affinization map. Its kernel $N$ is smooth with $h^0(\O_N)=1$, according to \cite{Demazure; Gabriel 1970}, Chapter III, \S3, 8.2.
It must be an extension of some abelian variety by a torus,
by \cite{Brion 2009}, Proposition 2.2. In our situation, the torus is trivial, so $N$ is an abelian variety.
The cokernel for $U\ra E^\aff$ is affine (loc.\ cit., Chapter III, \S3, Theorem 5.6). This cokernel is also the quotient of the abelian variety $A=E/U$.
It follows that $U\ra E^\aff$ is surjective, whence $E^\aff$ is unipotent, of some dimension $n\geq 1$.
The kernel for $U\ra E^\aff$ is affine and belongs to $N$, hence is proper, and therefore finite. In turn, $\dim(N)=\dim(A)$.
\qed

\section{The case of algebraic curves}
\mylabel{Algebraic curves}

Let $k$ be a ground field of characteristic $p\geq 0$. In this section,
$C$ denotes an \emph{algebraic curve}, that is, a scheme that is separated, of finite type,
equi-dimensional, and of dimension $d=1$.  
Write $C_1,\ldots,C_r$ for the irreducible components.   We regard each $C_i$
as the schematic images of the local Artin schemes $\Spec(\O_{C,\eta_i})$,
where  $\eta_i\in C$ denote the generic points.
Note that we do not assume that $C$ is proper, or reduced. However,   each $C_i$ is either proper or affine. 
We start by describing the affinization map for $C$:

\begin{proposition}
\mylabel{affinization curve}
The $k$-algebra $\Gamma(C,\O_C)$ is of finite type, the 
affinization map $f:C\ra C^\aff$ is projective with $f_*(\O_C)=\O_{C^\aff}$, and the exceptional locus
$\Exc(C/C^\aff)=\Supp(\Omega^1_{C/C^\aff})$ is the union of the irreducible components $C_i$ that are proper.
\end{proposition}

\proof
Choose some compactification $X=\bar{C}$. By definition, all embedded points lie inside the open set $C$,
so the local rings for the points at infinity are Cohen--Macaualay. Thus there is  some effective Cartier divisor $D\subset X$ whose support
is the closed set $X\smallsetminus C$. The ensuing short exact sequence $0\ra\shL^{\otimes n-1}\ra \shL^{\otimes n}\ra\shL_D^{\otimes n}\ra 0$
induces a long exact sequence
$$
H^0(X,\shL^{\otimes n})\lra H^0(X,\shL^{\otimes n}_D)\lra H^1(X, \shL^{\otimes n-1})\lra H^1(X, \shL^{\otimes n})\lra 0.
$$
We conclude that $h^1(\shL^{\otimes n})$ is decreasing in $n$, hence becomes constant for $n\gg 0$.
Then the map on the left must be surjective, and it follows that $\shL$ is semi-ample. Passing to a multiple
we may assume that $\shL$ is globally generated.
The homogeneous spectrum $Y=P(X,\shL)$ of the graded ring $R(X,\shL)=\bigoplus_{n\geq 0} \Gamma(X,\shL^{\otimes n})$
is a projective scheme and defines a morphism $f:X\ra Y$ with $\O_Y=f_*(\O_X)$. Moreover, $\shL$ is the preimage of the ample invertible sheaf $\O_Y(1)$.
The map contracts the irreducible components $X_i=\bar{C}_i$ that are disjoint from $D$; these  are exactly
the $C_i$ that are   proper. Let $s\in\Gamma(Y,\O_Y(1))$ be the global section whose preimage $f^*(s)\in\Gamma(X,\shL)$ vanishes precisely at $D$.
Then the non-zero locus $Y_s$ is an affine open set, with preimage $C=X\smallsetminus D=f^{-1}(Y_s)$.
Using $\O_Y=f_*(\O_X)$ we infer that $Y_s$ must be the affinization $C^\aff$.
\qed

\medskip
Suppose that the affine irreducible components $C_i$ are generically reduced.
Applying \cite{EGA II}, Corollary 7.4.11 to the normalization of $C_\red$, we see that there is a compactification $\bar{C}$
such that for each point at infinity $a\in \bar{C}$ the local ring $\O_{\bar{C},a}$ is a discrete valuation ring.
We call it the \emph{canonical compactification}. 

\begin{proposition}
\mylabel{canonical compactification}
Assumptions as above. Then the  canonical compactification $\bar{C}$ is an  initial object in the category $\Cpt(C)$
of all compactifications. Moreover, the formation of $\bar{C}$ commutes with separable ground field
extensions $k\subset k'$.
\end{proposition}

\proof
Let $C\subset X$ be any compactification. By assumption, the finite set $Z=X\smallsetminus C$ admits
an open neighborhood $U$ such that $U\smallsetminus Z$ is regular. It then follows that the canonical compactification
$\bar{C}$ arises from $X$ by normalization on the open set $U$, and making no change on the open set $X\smallsetminus Z$.
In particular, there is a morphism $\bar{C}\ra X$ between compactifications, so $\bar{C}$ yields an initial object in $\Cpt(C)$.

For the second assertion, write the locus at infinity as $\bar{C}\smallsetminus C=\{a_1,\ldots,a_r\}$,
and consider the residue fields $k_i=\kappa(a_i)$. Suppose the field extension $k\subset k'$ has the property
that the finite $k'$-algebra $k_i\otimes k'$ is regular, that is,  a product of fields. Then $\bar{C}\otimes k'$ remains regular over each $a_i\in \bar{C}$,
hence must coincide with the canonical compactification of $C\otimes k'$. This happens in particular if $k\subset k'$ is separable.
\qed

\medskip
We see that the canonical compactification $X_\lambda=\bar{C}$ is a final object in the opposite category $\Cpt(C)^\op$.
Choosing the singleton $L=\{\lambda\}$ as a cofinal index set, the proof for Theorem \ref{existence albanese map} immediately gives:

\begin{proposition}
\mylabel{albanese for curves}
Assumptions as in the previous proposition. Suppose $C^\aff$ is connected and reduced, and  $k=H^0(\bar{C},\O_{\bar{C}})$.
Then  the composition $C\subset \bar{C}\ra\Alb_{\bar{C}/k}$ is the Albanese map for the algebraic curve $C$.  
\end{proposition}

Note  that there is no final object in $\Cpt(C)^\op$ for more general curves $C$. For example,
the infinitesimal extensions $\PP^1\oplus\O_{\PP^1}(n)$ of the projective line, with $n$ arbitrary, can be seen as compactifications
of the non-reduced   curve $C= \AA^1_{k[\epsilon]}$, where $\epsilon^2=0$.

Let us unravel the condition that $C^\aff$ is connected and reduced.
Write  $\Gamma(C)$ for the \emph{dual graph} of the scheme $C$.
Recall that its  vertices correspond to the irreducible components $C_i$,
and two vertices are joined by an edge if $C_i\cap C_j$ is non-empty (compare for example \cite{Schroeer 2023}, discussion before Proposition 1.2). 
Let $C'\subset C$ be the union of the proper irreducible components, and 
$C''\subset C$ be the union of the affine irreducible components.
We regard these closed sets as closed subschemes, by declaring $C''$ as the schematic image of 
the morphism $\Spec(\prod\O_{C,\eta})\ra C$, where the product runs over the generic points $\eta\in C''$,
and likewise for $C'$.
These closed subschemes correspond to coherent sheaves of ideals $\shI'\subset\O_C$ and
$\shI''\subset\O_C$, respectively. Under the assumption that $C$ has no embedded components, we have  $\Supp(\shI'')=C'$, 
and may regard the sheaf of ideals $\shI''$  also as abelian sheaf on  $C'$, sitting in a short exact sequence
$0\ra \shI'\cap\shI''\ra \shI''\ra \shI''\O_{C'}\ra 0$, where the outer terms are $\O_{C'}$-modules.

\begin{proposition}
\mylabel{characterization affinization curve}
Notation as above. Suppose that the algebraic curve $C$ is not proper.
Then the affine hull $C^\aff$ is connected and reduced if and only if the following conditions hold:
\begin{enumerate}
\item The  dual graph $\Gamma(C)$ is connected.
\item The scheme $C$ has no embedded components.
\item The affine curve $C''$ is generically reduced.
\item The group  of global sections  $H^0(C',\shI'')$ vanishes.
\end{enumerate}
\end{proposition}

\proof
First, suppose that  that $C^\aff$ is connected and reduced. 
Then the noetherian scheme $C$ must be connected,   so the same holds for the dual graph $\Gamma(C)$.
Moreover, the structure sheaf $\O_C$ has no  global sections whose support is zero-dimensional.
Since $\dim(C)=1$ we infer that $C$ has no embedded components.
According to Proposition \ref{affinization curve}, the induced map $C''\ra C^\aff$ is a modification, and it follows
that $C''$ is generically reduced. The short exact sequence $0\ra\shI''\ra \O_C\ra\O_{C''}\ra 0$
gives an exact sequence
$$
0\lra H^0(C,\shI'')\lra H^0(C,\O_C)\lra H^0(C'',\O_{C''}),
$$
hence the term on the left vanishes. The latter can be written as  $H^0(C',\shI'')$, in light
of our chosen scheme structure on $C''$.

Conversely, suppose that conditions (i)--(iv) hold. The scheme $C$ is connected by (i).
The map $C\ra C^\aff$ is surjective, according to Proposition \ref{affinization curve}, so $C^\aff$ is connected as well.
In the above exact sequence, the term on the left vanishes by (iv), and the ring on the right is reduced by (iii).
Consequently  $H^0(C,\O_C)$ is reduced.
\qed

\medskip
Let $k\subset k'$ be a field extension, and set $C'=C\otimes k'$. Suppose that $C'$ is connected and reduced, and $k'=H^0(\bar{C}',\O_{\bar{C}'})$,
and consider the resulting comparison map 
\begin{equation}
\label{comparison curves}
c':\Alb_{C'/k'}\lra \Alb_{C/k}\otimes k'
\end{equation}
We shall see now that this may fail to be an  isomorphism.
As in \cite{Fanelli; Schroeer 2020b}, Section 2 we write $\Sing(\bar{C}/k)$ for the \emph{locus of non-smoothness}.
It carries a scheme structure, defined by the first Fitting ideal $\Fitt_1(\Omega^1_{\bar{C}/k})$. 
Note that it contains the \emph{singular locus} $\Sing(\bar{C})$, but over imperfect fields may be much larger.

\begin{theorem}
\mylabel{base change fails}
In the above setting, suppose that for the curve $C$ and the extension $k\subset k'$ the following holds:
\begin{enumerate}
\item The locus of non-smoothness   $\Sing(\bar{C}/k)$ is non-empty and contained in the locus at infinity $\bar{C}\smallsetminus C$.
\item The normalization  $X$ for the base-change $Y=\bar{C}\otimes k'$ is a smooth curve of genus $g\geq 1$.
\end{enumerate}
Then the comparison map \eqref{comparison curves} is not an isomorphism.
\end{theorem}

\proof
Our assumptions   ensure  that $X$ is  the canonical compactification of $C'$, and we have    a commutative diagram
$$
\begin{CD}
Y	@<\nu<<	X\\
@VgVV		@VVfV\\
\Alb_{Y/k'}	@<<c<	\Alb_{X/k'}.
\end{CD}
$$
By definition, $\Alb_{Y/k'}$ is the  base-change of $\Alb_{C/k}$, whereas $\Alb_{X/k}=\Alb_{C'/k'}$.  
The Albanese map $f:X\ra \Alb_{X/k'}$ is a closed embedding, in light of   (ii).

Seeking a contradiction, we assume that the comparison map is an isomorphism. 
Then the composition $g\circ \nu=c\circ f$ is affine, and the same holds for the normalization map $\nu:X\ra Y$.
Using the Leray--Serre spectral sequence and Serre's Criterion (\cite{EGA II}, Corollary 5.2.2), we infer that $g:Y\ra \Alb_{Y/k}$ is affine.
Since the composition 
$$
\O_{\Alb_{Y/k'}}\ra g_*(\O_Y)\ra (g\circ\nu)_*(\O_X)
$$
is surjective, it follows that $\nu:X\ra Y$ is a closed embedding. Hence this is an isomorphism,
consequently $\bar{C}$ is smooth, in contradiction to   (i).
\qed

\medskip
Let us discuss    explicit examples. Suppose the ground field $k$ is imperfect of characteristic $p>0$, 
and consider inside $\AA^2=\Spec k[x,y]$ the affine plane curve
$$
Z:\quad y^l = \prod_{i=1}^n (x^{q_i}-\lambda_i),
$$
where $l\geq 2$ is prime to  the characteristic,   and $q_i=p^{\nu_i}>1$ are powers of the characteristic, and $\lambda_i\in k$
are pairwise different scalars that are not $p$-powers in $k$, and  $n\geq 3$ is some   integer.
The right-hand side of the equation is  an inseparable   square-free polynomial of degree $q=q_1+\ldots+q_n$, 
and we see with Eisenstein's Criterion that $Z$ is integral. 
Note that such curves  where used by Totaro to construct  
\emph{pseudo-abelian varieties} that are not abelian varieties (\cite{Totaro 2013}, Example 3.1).

To simplify the exposition, we now also assume
that $l=q+1$. Then the projection $Z\ra \AA^1=\Spec k[x]$ immediately gives the canonical
compactification $h:\bar{Z}\ra\PP^1$, via the equation $y'^l=x'\cdot\prod(1-\lambda_ix'^{q_i})$ in the
new indeterminates    $x'=1/x$ and $y'=y/x$. One sees that $h^{-1}(\infty)$ contains only the  rational point given by $x'=y'=0$
and that $\bar{Z}$ is regular there, and hence also smooth (\cite{Fanelli; Schroeer 2020a}, Corollary 2.6).

The ideal   $\ideala\subset k[x,y]$ for the locus of non-smoothness $\Sing(Z/k)\subset \AA^2$ is generated by the defining equation and its partial derivatives.
One easily computes that the underlying closed set is given  by $y=0$ and $\prod(x^{q_i}-\lambda_i)=0$.
The equation $y=0$ defines an  effective Cartier divisor $D\subset Z$ whose coordinate ring  is the product of the
residue fields $k(\lambda_i^{1/q_i})$, and it follows that the affine plane curve $Z$ is regular.

Let $k\subset k'$ be an extension field that  contains the roots
$\lambda_i^{1/q_i}$,  and set $C'=C\otimes k'$.  At each singularity $a_i\in C'$, the complete local ring
becomes $R=k'[[u,v]]/(u^l-v^m)$, where $u=x$ and $v=y-\lambda_i^{1/q_i}$ and $m=q_i$.
Its normalization is $k[[t]]$, with normalization map   determined by  $u=t^m$ and $v=t^l$.
 We conclude that the normalization $X$ of the base-change $Y=\bar{Z}\otimes k'$
is smooth, of certain genus $g=h^1(\O_X)$. It comes with a branched covering $X\ra\PP^1_{k'}$ of degree $l$. In the Riemann--Hurwitz Formula  
$2g-2=l\cdot (-2) + \sum_{i=0}^n(l-1)$, the term on the right is
$$
l\cdot (-2) + \sum_{i=0}^n(l-1) = (n+1)(l-1)-2l \geq 4(l-1)-2l = 2l -4\geq 0,
$$
and it follows that $X$ is a curve of genus  $g\geq 1$.
Now consider the affine curve  $C=Z\smallsetminus \Sing(Z/k)$. Then $\bar{C}=\bar{Z}$,
and the two conditions in Theorem \ref{base change fails} are satisfied. We see that the comparison map \eqref{comparison curves} is not an isomorphism.

\section{The case of algebraic groups}
\mylabel{Algebraic groups}

Let $k$ be a ground field of characteristic $p\geq 0$.  Throughout this section, $G$
denotes a group scheme of finite type, which are also called \emph{algebraic groups} in the literature.
First note that $G$ is connected and reduced if and only if the respective properties hold for the
\emph{affinization} $G^\aff=\Spec\Gamma(G,\O_G)$. Throughout, we will assume that these equivalent conditions hold.

Since the neutral element $e\in G$ is a rational point, the ground field $k$ coincides with the essential field of constants
for $G$. By Theorem \ref{existence albanese map}, there is an  Albanese map $f:G\ra\Alb_{G/k}$.
The Albanese variety is a para-abelian variety endowed with a rational point $f(e)$.
The latter becomes the zero element $0\in\Alb_{G/k}$ for a unique group law, according to  \cite{Laurent; Schroeer 2021}, Proposition 4.3.
In what follows we regard $\Alb_{G/k}$ as an \emph{abelian variety}. 
 
The goal of this section is to analyze how the  various group laws
interact with the Albanese map.  We shall see  that   $f$ does not necessarily respect the group laws.
However, the    following key fact (\cite{Brion 2017}, Proposition 4.1.4) ensures that it does  so on many 
subgroup schemes:

\begin{lemma}
\mylabel{restriction respects group law}
The restriction of the Albanese map $f:G\ra \Alb_{G/k}$ to any smooth connected subgroup scheme $H\subset G$ respects the group laws.
\end{lemma}

If $k$ is perfect, then the reduced group scheme $G$ itself is smooth, so the Albanese map is the universal homomorphism 
to an abelian variety. 

The following terminology, which applies to any closed subscheme  $E\subset G$  containing the origin $e\in G$, will be useful:
We say that $f|E$ is \emph{trivial} if it factors over the zero element $0\in\Alb_{G/k}$, viewed as a reduced closed subscheme.
More generally, we say that $f|E$ is \emph{set-theoretically trivial} if $f(g)=0$ for every  point $g\in E$.
In this case,  the \emph{schematic image} $Z\subset \Alb_{G/k}$ is some closed subscheme supported by the zero element.
In turn, $Z=\Spec(R)$ for some finite local $k$-algebra $R$ with residue field $R/\maxid_R=k$.

\begin{proposition}
\mylabel{restriction is trivial}
The restriction $f|H$  to any reduced connected affine subgroup scheme $H\subset G$   is   trivial.
\end{proposition}

\proof
We start with the case that $f|H$ respects the group law.
According to \cite{Demazure; Gabriel 1970}, Chapter II, Proposition 5.1 the set-theoretical image of $f|H$  is a closed 
set. It must be connected, because $H$ is connected. Write $E$ for this closed set, endowed with the reduced scheme structure. 
Since $H$ is reduced,   $f:H\ra\Alb_{G/k}$ factors over $E$.
By loc.\ cit.\ $E$ is a subgroup scheme, and the homomorphism $f:H\ra E$ is faithfully flat.
In turn, we have $E=H/H'$ where $H'=\Kernel(f|H)$, and the quotient must be affine, according to \cite{Demazure; Gabriel 1970}, Chapter III, Theorem
in 7.2. Since the Albanese variety is proper, the same holds for $E$.
Being proper and affine, the group scheme $E$ must be finite.
Since it is connected it must be a singleton, thus  $f|H$ is trivial. 

In light of Lemma \ref{restriction respects group law}, our assertion holds if $H$ is smooth.
Suppose now that $H$ is not smooth. Then we are in characteristic $p>0$.
Consider the \emph{Frobenius kernels} $H[F^n]$ for the iterated relative Frobenius maps
$F^n:H\ra H^{(p^n)}$. According to \cite{Demazure; Gabriel 1970}, Chapter III, Lemma 6.10  the quotient $H/H[F^n]$ is smooth
for sufficiently large $n\geq 0$. Set $A=\Alb_{G/k}$ and consider the induced morphism $f_n:G^{(p^n)}\ra A^{(p^n)}$.
Its restriction to the subgroup scheme $H/H[F^n]$ inside $ H^{(p^n)}\subset G^{(p^n)}$ is  trivial, by the preceding paragraph,
and we infer that $f|H$ is at least set-theoretically trivial.
Let $Z=\Spec(R)$ be its schematic image. Then $R$ is a finite local $k$-algebra with residue field $R/\maxid_R=k$,
and the canonical map $R\ra\Gamma(H,\O_H)$ is injective. Since $H$ is reduced, the same holds for $R$,
and thus $R=k$.
\qed

\medskip
Note that the above arguments also give the fact  that each subscheme $H\subset G$ that inherits a group law must be a closed subscheme.
Each such $H\subset G$ acts via translations  on $G$, from both sides.
If $f|H$ respects the group law,  $H$ likewise acts on $\Alb_{G/k}$   via translations, also from both sides.

\begin{proposition}
\mylabel{albanese map is equivariant}
Let $H\subset G$ be a smooth connected subgroup scheme.  
Then the Albanese map $f:G\ra\Alb_{G/k}$ is equivariant with respect to the $H$-actions.
\end{proposition}

\proof
We have to verify that the graph $\Gamma_f\subset G\times\Alb_{G/k}$
is $H$-stable, where $H$ acts diagonally. In light of Theorem \ref{base change}, it suffices to treat the case that $k$ is separably closed.
Let $H'\subset H$ be the stabilizer of the graph.
Each element $\sigma\in H(k)$ defines a translation automorphism $\sigma:G\ra G$.
From the universal property of Albanese maps, it induces a compatible automorphism  
of the Albanese variety; it follows that   the inclusion $H'(k)\subset H(k)$ is an equality.
Since $H$ is smooth, the set $H(k)$ is dense, and thus the inclusion of the closed set  $H'$ is an equality.
\qed

\medskip
Note that for reduced but non-smooth $H$, the regular locus $\Reg(H)$ is open and dense, but contains
no rational point (\cite{Fanelli; Schroeer 2020a}, Corollary 2.6), and it can  easily happen that the group  $H(k)$
is trivial (\cite{Schroeer; Tziolas 2023}, Section 8).

\begin{corollary}
\mylabel{albanese map factors}
Let $H\subset G$ be a smooth connected subgroup scheme that is also  affine. Then the Albanese map $f:G\ra \Alb_{G/k}$ uniquely factors over $G/H$.
\end{corollary}

\proof
By Proposition \ref{restriction is trivial}, the homomorphism $f|H$ is  trivial, hence the $H$-action on the Albanese variety is trivial.
Passing to quotients we get  $G/H\ra\Alb_{G/k}$. Such a factorization is unique because the
projection $G\ra G/H$ is an epimorphism.
\qed

\medskip
Since   affinization  of schemes that are separated and of finite type commutes with products,   $G^\aff$
inherits a group law, and the affinization map  $p:G\ra G^\aff$ is a homomorphism.
Let $N$ be its kernel and write $i:N\ra G$ for the inclusion map.
By the  Affinization Theorem (\cite{Demazure; Gabriel 1970}, Chapter III, 8.2),
we get a central extension
\begin{equation}
\label{central extension}
0\lra N\stackrel{i}{\lra} G\stackrel{p}{\lra} G^\aff\lra 1.
\end{equation}
Furthermore, $N$ is smooth and  connected, with $h^0(\O_N)=1$.
Such group schemes are called \emph{anti-affine}. Their structure was analyzed by Brion \cite{Brion 2009}.
Note that on $G$ and $G^\aff$, the group laws are written in a multiplicative way, whereas for $N$ we choose additive notation.
Let us say  that the group scheme $G$ \emph{weakly splits}   if there is a  morphism of schemes $s:G^\aff\ra G$ with
$p\circ s=\id_{G^\aff}$. Then
\begin{equation}
\label{identification weakly split}
N\times G^\aff\lra G,\quad (x,y)\longmapsto i(x)\cdot s(y)
\end{equation}
is an isomorphism of schemes. 
If the section $s$ additionally respects  the group laws, 
we say that $G$ \emph{strongly splits}; then the above is actually an isomorphism of group schemes.

We now can formulate our main result on the Albanese variety of  group schemes $G$ of finite type that are reduced and connected.
It unravels how the Albanese map $f:G\ra\Alb_{G/k}$ is related to the
central extension $0\ra N\ra G\ra G^\aff\ra 0$. Note that by
Proposition \ref{restriction is trivial}, the restriction $f|N$ respects the group law, and we set $N'=\Kernel(f|N)$.

\begin{theorem}
\mylabel{albanese for group schemes}
In the above situation, the kernel $N'\subset N$   inside the anti-affine group scheme
$N$ is the smallest subgroup scheme such that $N/N'$ is proper and   $G/N'$ weakly splits.
Moreover,  for any section $s:G^\aff\ra G/N'$, the composition
$$
G\stackrel{\can}{\lra} G/N'\stackrel{(i,s)^{-1}}{\lra} N/N'\times G^\aff\stackrel{\pr_1}{\lra} N/N'  
$$
is an Albanese map for $G$.
\end{theorem}

\proof
Let $B$ be the cokernel of $f|N$, which is an abelian variety.
According to Proposition \ref{albanese map is equivariant},
the Albanese map $f:G\ra\Alb_{G/k}$ induces a morphism $G^\aff=G/N\ra B$. The latter is trivial, by Proposition \ref{restriction is trivial},
and the universal property of $f$ reveals that $B=0$. Thus $f|N$ is surjective, and we get $\Alb_{G/k}=N/N'$.
In particular, $N/N'$ is proper.

Set $A=\Alb_{G/k}$. The Albanese map $f:G\ra A$    factors over $G/N'$.
The induced morphism $g:G/N'\ra A$ is equivariant with respect to the actions of $N/N'$
and the restriction of $g$ to $N/N'$ is an isomorphism of abelian varieties. In turn, $r=(g|N/N')^{-1}\circ g$ is
an \emph{equivariant  retraction} for the inclusion $j: N/N'\ra G/N'$. It follows
that the composition
$$
G^\aff\lra \{e\}\times G^\aff\lra N/N'\times G^\aff \stackrel{(r,p)^{-1}}{\lra} G
$$
is a section for the projection $G/N'\ra G^\aff$. Thus $G/N'$ weakly splits.
 
Next, we describe the Albanese map in  the special case that $G$ is weakly split. Since Albanese maps depend only on the underlying scheme,
it suffices to treat the case that $G$ is strongly split.  Let $N_1\subset N$ be the maximal smooth connected affine subgroup scheme.
Using that the Albanese map factors over $N/N_1$, we reduce to the case $N_1=0$. Then $N$ is an abelian variety,
(this depends on Brion's result \cite{Brion 2009}, see also  \cite{Laurent; Schroeer 2021}, Proposition 7.1).
We now have to verify that  $\pr_1:N\times G^\aff\ra N$ is an Albanese map. Consider the commutative diagram
$$
\begin{tikzcd} 
N\ar[r,"i"]\ar[dr]	& N\times G^\aff\ar[dr,"\pr_1"]\ar[d,"f"]\\
		& \Alb_{G/k}\ar[r,"g"']			 &   N,
\end{tikzcd}
$$
where $g$ comes from the universal property of the Albanese map. The left diagonal arrow $f\circ i$
is surjective, by the preceding paragraph, and has trivial kernel, because $g\circ f\circ i=\id_N$.
It follows that $g$ is an isomorphism, whence $\pr_1$ is an Albanese map.

We come back to the general case. Let $H\subset N$ be any subgroup scheme with  $G/H$ weakly split and $N/H$ proper.
It only remains to  check   that $N'\subset H$. Choose a section $s$ for the projection $G/H\ra G^\aff$ and consider the commutative diagram
$$
\begin{tikzcd} 
N\ar[r,"i"]\ar[d,"\can"']	& G\ar[r,"\can"]\ar[d,"f"]	& G/H\ar[r,"(\text{$i,s$})^{-1}" ]		& N/H\times G^\aff\ar[d,"\pr_1"]\\
N/N'\ar[r]	& \Alb_{G/k}\ar[rr,"g"']	& 				& N/H,
\end{tikzcd}
$$
where $g$ comes from the universal property of the Albanese map $f$. We see that the composite map $N\ra N/H$ is equivariant with respect to the action of $N$,
and factors over $N/N'$, which ensures $N'\subset H$.
\qed

\medskip
Suppose $G$ is weakly split, and choose a section $s$ for the projection $p$.
Then the identification   $i\cdot s: N\times G^\aff\ra G$ from \eqref{identification weakly split}
is an isomorphism of schemes, and the group law on $G$ arises from the product group law by a modification
with a Hochschild two-cocycle $\varphi:(G^\aff)^2\ra N$. In fact, the isomorphism classes of central extensions $0\ra N\ra E\ra G^\aff\ra 0$ where the projection
$E\ra G^\aff$ admits a section correspond to 
classes in the Hochschild cohomology group $H^2_0(G^\aff, N)$, where $N$ is viewed as a module over $G^\aff$ with trivial action, see the discussion in
\cite{Demazure; Gabriel 1970}, Chapter II, \S3. This yields the following:
 
\begin{corollary}
\mylabel{albanese respect group law}
Suppose our group scheme $G$ is weakly split and that $N$ is proper. 
Then the Albanese map $f:G\ra \Alb_{G/k}$ respects the group law if and only if   $G$ strongly splits.
\end{corollary}

\proof
We may assume that $G=N\times G^\aff$ as schemes, and choose the two-cocycle $\varphi:(G^\aff)^2\ra N$ 
so that the inclusion $i:N\ra G$ is given by $x\mapsto (x,e)$.
The projection $\pr_1:G\ra N$ is   an Albanese map, by the theorem.

Suppose now that $\pr_1:G\ra N$ respects the group law.
Then the identity morphism $G\ra N\times G^\aff$ respects the group laws, hence $G$ is strongly split.
Conversely, if $G$ is strongly split, then $\pr_1$ respects the group law.
\qed

\medskip
We now construct examples where the Albanese map does not respect the group laws.
Recall that a group scheme $U$ is \emph{unipotent} if it is of finite type,
and the base-change $U\otimes k^\alg$ admits  a   composition series whose quotients can be embedded into the additive group $\GG_{a, k^\alg}$.
We refer to  \cite{SGA 3b}, Expos\'e XVII for a  comprehensive treatment.
Note that if $U$ is connected,  at least one composition series is already defined over $k$ (loc.\ cit.\ Theorem 3.5).

\begin{proposition}
\mylabel{albanese disrespects group law}
Let $U$ be a reduced connected unipotent group scheme, and $N$ be an abelian variety.
Suppose we are in characteristic $p>0$, and that there is an epimorphism $U\otimes k^\alg\ra \alpha_{p,k^\alg}$
and a monomorphism $\alpha_{p,k^\alg}\ra N\otimes k^\alg$. Then the following holds:
\begin{enumerate}
\item The Hochschild cohomology group $H^2_0(U,N)$ is non-zero.
\item For all $\gamma\in H^2_0(U,N)$   the resulting extension $G$ is reduced and connected.
\item If $\gamma\neq 0$, the Albanese map $G\ra\Alb_{G/k}$ does not respect the group law.
\end{enumerate}
\end{proposition}

\proof
To see (ii), observe that in the extension $0\ra N\ra G\ra U\ra 1$ the projection $G\ra U$
has the structure of  a principal $N$-bundle. By our assumptions, all fibers are smooth and connected,
and the base is reduced and connected. So the latter also holds for the total space $G$. This gives (ii).
Assertion (iii) is a direct consequence of Corollary \ref{albanese respect group law}.

It remains to verify (i).
For this we first recall the computation of  the cocycle group  $Z^2(\GG_a,\GG_a)$ for Hochschild cohomology.
The cochains are just morphisms $\GG_a^2\ra\GG_a$. These correspond to homomorphisms
$k[Z]\ra k[X]\otimes k[Y]$ of $k$-algebras. The latter are determined by the images
of $Z$, hence are given by a polynomial in $X$ and $Y$.
By the computation in \cite{Demazure; Gabriel 1970}, Chapter II, \S3, Subsection 4.6
the polynomials
$$
Q(X,Y)=\sum_{i=1}^{p-1}\binom{p}{i} X^iY^{p-i}\quadand XY^{p^r} \quad (r\geq 1)
$$
satisfy the cocycle condition. We remark in passing that these freely generate a complement
for $B^2\subset Z^2$, viewed as modules over the skew polynomial ring $k[F]$, where $F\lambda=\lambda^pF$ holds.
The polynomial $Q(X,Y)$ induces a non-zero homomorphism
$$
k[Z]/(Z^p)\lra k[X]/(X^p)\otimes k[Y]/(Y^p),
$$
giving a non-trivial element in $Z^2(\alpha_p,\alpha_p)$.
Using  the epimorphism $U\otimes k^\alg\ra \alpha_{p,k^\alg}$ over $k^\alg$ and applying descent,
we conclude that $Z^2(U,\alpha_p)$ is non-zero. 

Now let us examine the group of one-cochains $C^1(U,\alpha_p)$, whose elements are morphism of schemes $h:U\ra\alpha_p$.
Since $U$ is reduced, any such morphism factors over  $e\in \alpha_p$, viewed as a reduced closed subscheme.
Thus the coboundary operator $C^1(U,\alpha_p)\ra Z^2(U,\alpha_p)$ is the zero map, and we conclude  $H^2_0(U,\alpha_p)\neq 0$.

Now consider the abelian variety $N$. Its Lie algebra $\Lie(N)$ has trivial brackets, and comes with an additional structure
given by the $p$-map $x\mapsto x^{[p]}$. By the Demazure--Gabriel Correspondence (\cite{Demazure; Gabriel 1970}, Chapter II, \S7, Theorem 3.5)
the homomorphisms
$\alpha_p\ra N$ correspond to the vectors $x\in \Lie(N)$ with $x^{[p]}=0$. Since the brackets are trivial, the $p$-map is additive.
 Our assumption on $N\otimes k^\alg$ ensures that
there is a non-zero homomorphism $\alpha_p\ra N$.
In turn, we get an inclusion $Z^2(U,\alpha_p)\ra Z^2(U,N)$, so these groups are non-zero.
Again we consider $C^1(U,N)$. Fix an element $h:U\ra N$. By Proposition \ref{restriction is trivial},
this morphisms factors over a rational point $a\in N$. In turn, the coboundary map $C^1(U,N)\ra Z^2(U,N)$ is trivial,
and we conclude again that $H^2_0(U,N)\neq 0$.
\qed

\medskip
It is   easy to make this concrete,  over any imperfect field $k$ of characteristic $p>0$:
Choose an element $t\in k$ that is not a $p$-power, and let $U$ be the kernel of the homomorphism
$\GG_a^2\ra \GG_a$ given by $(x,y)\mapsto x^p-ty^p$.
Then $U$ is integral,  $\pr_1|U$ is surjective, and the kernel is isomorphic to $\alpha_p$.
Over the field extension $k'=k(t^{1/p})$,   the section given by $z\mapsto (z,t^{1/p}z)$ 
defines the desired splitting $U\otimes k'\simeq (\GG_a\times \alpha_p)\otimes k'$.
Furthermore, there is a supersingular elliptic curve $N$ over $k$ (well-known, see for example \cite{Roessler; Schroeer 2022}, Lemma 3.1),
which indeed contains a copy of $\alpha_p$.


\end{document}